\documentclass[draft=false,bibliography=totoc]{scrartcl}

\usepackage{amsmath}
\usepackage{amssymb}
\usepackage{amsthm}
\usepackage{cite}
\usepackage[T1]{fontenc}
\usepackage{hyperref}
\usepackage{lmodern}
\usepackage{mathscinet}
\usepackage{mathtools}
\usepackage{microtype}
\usepackage[automark]{scrpage2}

\usepackage{string_arxiv}

 \newtheorem{thm}{Theorem}[section]

\theoremstyle{definition}
 \newtheorem{defn}[thm]{Definition}
 
\theoremstyle{remark}
 \newtheorem{rem}[thm]{Remark}

\mathtoolsset{mathic=true}

\numberwithin{equation}{section}
\title{On resolvent matrix, Dyukarev-Stieltjes parameters and orthogonal matrix polynomials via $\ra$-Stieltjes transformed sequences}
\author{Abdon Eddy Choque-Rivero\footnote{A.~E.~Choque-Rivero is supported by SNI--CONACYT and CIC--UMSNH, M\'exico.} \and Conrad M\"adler}

\begin{document}
\maketitle

\begin{abstract}
 By using \tKt{ed} sequences and Dyukarev-Stieltjes parameters we obtain a new representation of the resolvent matrix corresponding to the truncated matricial Stieltjes moment problem. Explicit relations between orthogonal matrix polynomials and matrix polynomials of the second kind constructed from consecutive \tKt{ed} sequences are obtained. Additionally, a \tnnH{} measure for which the matrix polynomials of the second kind are the orthogonal matrix polynomials is found.
\end{abstract}

\subsection*{Keywords:}
 Resolvent matrix, orthogonal matrix polynomials, Dyukarev-Stieltjes parameters, \tKt{ed} sequences.

\section{Introduction}
 This paper is a continuation of work done in the papers~\zitas{MR3324594,MR3327132,MR2735313,MR3014201,MR3133464}, where two truncated matricial power moment problems on semi-infinite intervals made up one of the main topics. The starting point of studying such problems was the famous two part memoir of \Stieltjes{}~\zitas{MR1508159,MR1508160} where the author's investigations on questions for special continued fractions led him to the power moment problem on the interval \(\ra\). A complete theory of the treatment of power moment problems on semi-infinite intervals in the scalar case was developed by M.~G.~Kre{\u\i}n in collaboration with A.~A.~Nudelman (see~\zitaa{MR0044591}{\cSect{10}},~\zita{MR0233157},~\zitaa{MR0458081}{\cchap{V}}). For a modern operator-theoretical treatment of the power moment problems named after Hamburger and Stieltjes and its interrelations, we refer the reader to Simon~\zita{MR1627806}.

 The matrix version of the classical \Stieltjes{} moment problem was studied in Adamyan/Tkachenko~\zitas{MR2155645,MR2215856}, And\^o~\zita{MR0290157}, Bolotnikov~\zitas{MR975671,MR1362524,MR1433234}, Bolotnikov/Sakhnovich~\zita{MR1722780}, Chen/Hu~\zita{MR1807884}, Chen/Li~\zita{MR1670527}, Dyukarev~\zitas{Dyu81,MR686076}, Dyukarev/Katsnel{\('\)}son~\zitas{MR645305,MR752057}, Hu/Chen~\zita{MR2038751}.


 The central research object of the present work is the resolvent matrix (RM) $\Uu{m}$ of the truncated matricial Stieltjes matrix moment (TSMM) problem. The importance of the knowledge of the RM $\Uu{m}$ is explained by the fact that the matrix $\Uu{m}$ generates the solution set of the TSMM problem via a linear fractional transformation. The multiplicative decomposition of $\Uu{m}$ in simplest factors containing Dyukarev-Stieltjes (DS) parameters $\DSpm{j}$ and $\DSpl{j}$~\zita{MR3324594} allowed us to attain interesting interrelations between the orthogonal polynomials $P_{k,j}$, their second kind polynomials $Q_{k,j}$ and the DS-parameters $\DSpm{j}$ and $\DSpl{j}$, as well as the Schur complements $\widehat H_{k,j}=L_{k,j}$; see~\zita{MR3324594}.

 In the present work, inspired by Dyukarev's multiplicative decomposition of \(\Uu{m}\) (see \rprop{P1033}), the \tKt{ed} sequences \zitas{114arxiv,141arxiv}, and the representation of the RM in terms of the polynomials $P_{k,j}$ and $Q_{k,j}$~\zita{MR3324594}, a factorization of the RM $\Uu{m}$ is obtained. This representation is constructed through a sequence of \taKt{\ell}ed sequences and corresponding polynomials $P_{k,j}^{(\ell)}$ and $Q_{k,j}^{(\ell)}$. An important consequence of such representation is the fact that a \tnnH{} measure for which the polynomials of the second kind  $Q_{k,j}^{(1)}$ are orthogonal is explained. By employing the interrelations between the orthogonal matrix polynomials and the Hurwitz type matrix polynomials (see~\cite{MR3327132}) new identities involving $P_{k,j}$, $Q_{k,j}$ and $P_{k,j}^{(1)}$, $Q_{k,j}^{(1)}$ are attained; see \rtheo{Tnew1}. The scalar version of the mentioned interrelations were studied in~\cite{CR16}.


 The starting point of our considerations in this paper is the \tDSp{} \([\seq{\DSpl{k}}{k}{0}{\infty},\seq{\DSpm{k}}{k}{0}{\infty}]\) of a matricial moment sequence corresponding to a completely \tnd{} \tnnH{} measure on the right half-axis \(\ra\), having moments up to any order. Yu.~M.~Dyukarev~\zita{MR2053150} introduced these parameters in connection with a multiplicative decomposition of his resolvent matrix for the truncated matricial Stieltjes moment problem into elementary factors to characterize indeterminacy. In the scalar case the \tDSp{} coincides with the classical parameters \([\seq{\spl{k}}{k}{0}{\infty},\seq{\spm{k}}{k}{0}{\infty}]\) used by \Stieltjes{}~\zitas{MR1508159,MR1508160} to formulate his indeterminacy criterion. M.~G.~Kre{\u\i}n gave a mechanical interpretation for \Stieltjes{'} investigations on continued fractions (see Gantmacher/Kre{\u\i}n~\zitaa{MR0114338}{\canha{2}} or Akhiezer~\zitaa{MR0184042}{Appendix}) by a weightless thread carrying point masses \(m_k\) with intermediate distances \(l_k\). In~\zita{MR3014201} another parametrization of moment sequences related to a semi-infinite interval \([\alpha,\infp)\) was introduced, the so-called \tSp{} \(\seq{\Spu{j}}{j}{0}{\infi}\). This parametrization is strongly connected to a Schur-type algorithm considered in~\zitas{MR2038751,MR1807884,114arxiv,141arxiv} for solving the truncated matricial Stieltjes moment problem step-by-step by reducing the number of given data. In one step the sequence of prescribed moments is transformed into a shorter sequence, the so-called \tKt{}, eliminating one moment. This procedure is equivalent to dropping the first Stieltjes parameter \(\Spu{0}\). In \rtheo{R0921} we show that this transformation is also essentially equivalent to dropping the Dyukarev-Stieltjes parameter \(\DSpm{0}\) and interchanging the roles of \(\DSpl{k}\) and \(\DSpm{k}\),
  which is especially interesting against the background of M.~G.~Kre{\u\i}n's mechanical interpretation. By dividing the elementary factors of Yu.~M.~Dyukarev's above mentioned multiplicative representation of his resolvent matrix into two groups, in \rtheo{Taaa007} we give a factorization of this resolvent matrix, which corresponds to a splitting up of the original problem into two smaller moment problems associated with the first part of the original sequence of prescribed moments and a second part obtained by repeated application of the \tKt{ation}. Comparing blocks in this formula, in \rtheo{T0820} we can represent the orthogonal matrix polynomials and second kind matrix polynomials with respect to the \tKt{ed} moment sequence in terms of the polynomials corresponding to the original moment sequence. In particular, we state orthogonality relations for the matrix polynomials of the second kind in \rprop{orthg-001}.

 In order to formulate the moment problems we are going to study, we first review some notation. Let \(\C\)\index{c@\(\C\)}, \(\R\)\index{r@\(\R\)}, \(\NO\)\index{n@\(\NO\)}, and \(\N\)\index{n@\(\N\)} be the set of all complex numbers, the set of all real numbers, the set of all \tnn{} integers, and the set of all positive integers, respectively. Throughout this paper, let \(p,q\in\N\)\index{p@\(p\)}\index{q@\(q\)}. For all \(\alpha,\beta\in\R\cup\set{-\infty,\infp}\), let \(\mn{\alpha}{\beta}\)\index{z@\(\mn{\alpha}{\beta}\)} be the set of all integers \(k\) for which \(\alpha\leq k\leq\beta\) holds. If \(\cX\) is a nonempty set, then \(\cX^\pxq\)\index{\(\cX^\pxq\)} stands for the set of all \tpqa{matrices}, each entry of which belongs to \(\cX\), and \(\cX^p\)\index{\(\cX^p\)} is short for \(\cX^\xx{p}{1}\). If \((\Omega,\gA)\) is a measurable space, then each countably additive mapping whose domain is \(\gA\) and whose values belong to the set \(\Cggq\)\index{c@\(\Cggq\)} of all \tnn{} \tH{} complex \tqqa{matrices} is called a \tnnH{} \tqqa{measure} on \((\Omega,\gA)\). Denote by \symba{\Cgq}{c} the set of all \tpH{} complex \tqqa{matrices}.

 Let \(\Bra\)\index{b@\(\Bra\)} be the \(\sigma\)\nobreakdash-algebra of all Borel subsets of \(\ra\), let \(\Mggqa{\ra}\)\index{m@\(\Mggqa{\ra}\)} be the set of all \tnn{} \tH{} \tqqa{measures} on \((\ra,\Bra)\) and, for all \(\kappa\in\NOinf \), let \(\Mgguqa{\kappa}{\ra}\)\index{m@\(\Mgguqa{\kappa}{\ra}\)} be the set of all \(\sigma\in\Mggqa{\ra}\) such that the integral
 \bgl{so}
  \suo{j}{\sigma}
  \defg\int_{\ra}t^j\sigma(\dif t)
 \eg
 \index{s@\(\suo{j}{\sigma}\)}exists for all \(j\in\mn{0}{\kappa}\).
 Two matricial power moment problems lie in the background of our considerations. The first one is the following:
 \begin{description}
  \item[\mproblem{\ra}{\kappa}{=}]\index{m@\mproblem{\ra}{\kappa}{=}} Let \(\kappa\in\NOinf \) and let \(\seqska \) be a sequence of complex \tqqa{matrices}. Describe the set \(\Mggqaag{\ra}{\seqska }\)\index{m@\(\Mggqaag{\ra}{\seqska }\)} of all \(\sigma\in\Mgguqa{\kappa}{\ra}\) for which \(\suo{j}{\sigma}=\su{j}\) is fulfilled for all \(j\in\mn{0}{\kappa}\).
 \end{description}
 The second matricial moment problem under consideration is a truncated one with an additional inequality condition for the last prescribed moment:
 \begin{description}
  \item[\mproblem{\ra}{m}{\leq}]\index{m@\mproblem{\ra}{m}{\leq}} Let \(m\in\NO\) and let \(\seq{\su{j}}{j}{0}{m}\) be a sequence of complex \tqqa{matrices}. Describe the set \(\Mggqaakg{\ra}{\seq{\su{j}}{j}{0}{m}}\)\index{m@\(\Mggqaakg{\ra}{\seq{\su{j}}{j}{0}{m}}\)} of all \(\sigma\in\Mgguqa{m}{\ra}\) for which \(\su{m}-\suo{m}{\sigma}\) is \tnn{} \tH{} and, in the case \(m>0\), moreover \(\suo{j}{\sigma}=\su{j}\) is fulfilled for all \(j\in\mn{0}{m-1}\).
 \end{description}

In order to give a better motivation for our considerations in this paper,
 we are going to recall the characterizations of solvability of the above mentioned moment problems, which were obtained in~\zita{MR2735313}.
  This requires some preparations.

For all \(n\in\NO\), let \(\Kggqu{2n+1}\)\index{k@\(\Kggqu{2n+1}\)} be the set of all sequences \(\seq{\su{j}}{j}{0}{2n+1}\) of complex \tqqa{matrices}, such that the block \Hankel{} matrices
\begin{align}\label{HK}
 \Hu{n}&\defg\matauuo{\su{j+k}}{j,k}{0}{n}&
 &\text{and}&
 \Ku{n}&\defg\matauuo{\su{j+k+1}}{j,k}{0}{n}
\end{align}
\index{h@\(\Hu{n}\)}\index{k@\(\Ku{n}\)}are both \tnnH{}. For all \(n\in\NO\), let \(\Kggqu{2n}\)\index{k@\(\Kggqu{2n}\)} be the set of all sequences \(\seq{\su{j}}{j}{0}{2n}\) of complex \tqqa{matrices}, such that \(\Hu{n}\) is \tnnH{} and, in the case \(n\geq1\), furthermore \(\Ku{n-1}\) is \tnnH{}. For all \(m\in\NO\), let \(\Kggequ{m}\)\index{k@\(\Kggequ{m}\)} be the set of all sequences \(\seq{\su{j}}{j}{0}{m}\) of complex \tqqa{matrices} for which a complex \tqqa{matrix} \(\su{m+1}\) exists such that \(\seq{\su{j}}{j}{0}{m+1}\) belongs to \(\Kggqu{m+1}\). Let \(\Kggqu{\infty}\)\index{k@\(\Kggqinf\)} be the set of all sequences \(\seqsinf \) of complex \tqqa{matrices} such that \(\seq{\su{j}}{j}{0}{m}\) belongs to \(\Kggqu{m}\) for all \(m\in\NO\), and let \(\Kggequ{\infty}\defg\Kggqu{\infty}\)\index{k@\(\Kggeqinf\)}. For all \(\kappa\in\NOinf \), we call a sequence \(\seqska \) \noti{\tSnnd}{Stieltjes!\tnn{} definite} (resp.\ \noti{\tSnnde}{Stieltjes!\tnn{} definite!extendable}) if it belongs to \(\Kggqu{\kappa}\) (resp.\ to \(\Kggequ{\kappa}\)). Observe that these notions coincide with the right-sided version in~\zitaa{MR3014201}{\cdefnp{1.3}{213}} for \(\alpha=0\).

 Using the sets of matrix sequences above, we are able to formulate the solvability criterions of the problems~\mproblem{\ra}{m}{=} and~\mproblem{\ra}{m}{\leq}, which were obtained in~\zita{MR2735313} for intervals \([\alpha,\infp)\) with arbitrary \(\alpha\in\R\):
\btheol{T1121}
 Let \(\kappa\in\NOinf \) and let \(\seqska \) be a sequence of complex \tqqa{matrices}.
  Then \(\Mggqaag{\ra}{\seqska }\neq\emptyset\) if and only if \(\seqska \in\Kggequ{\kappa}\).
\etheo
 In the case \(\kappa\in\NO\), \rtheo{T1121} is a special case of~\zitaa{MR2735313}{\cthmp{1.3}{909}}. If \(\kappa=\infi\), the asserted equivalence can be proved using the equation \(\Mggqaag{\ra }{\seqsinf }=\bigcap_{m=0}^\infty\Mggqaag{\ra }{\seq{\su{j}}{j}{0}{m}}\) and the matricial version of the Helly-Prohorov theorem (see~\zitaa{MR975253}{\cSatzp{9}{388}}). We omit the details of the proof, the essential idea of which is originated in~\zitaa{MR0184042}{proof of \cthmp{2.1.1}{30}}.

\bthmnl{see~\zitaa{MR2735313}{\cthmp{1.4}{909}}}{T1122}
 Let \(m\in\NO\) and let \(\seq{\su{j}}{j}{0}{m}\) be a sequence of complex \tqqa{matrices}. Then \(\Mggqaakg{\ra}{\seq{\su{j}}{j}{0}{m}}\neq\emptyset\) if and only if \(\seq{\su{j}}{j}{0}{m}\in\Kggqu{m}\).
\end{thm}

 The importance of \rtheoss{T1121}{T1122} led us in~\zitas{MR2735313,MR3014201,MR3133464} to a closer look at the properties of sequences of complex \tqqa{matrices}, which are \tSnnd{} or \tSnnde{}. Guided by our former investigations on \tHnnd{} sequences and \tHnnde{} sequences, which were done in~\zita{MR2570113}, in~\zitaa{MR2735313}{\cSect{4}} and in~\zita{MR3014201} we started a thorough study of the structure of \tSnnd{} sequences and \tSnnde{} sequences.


 This paper is organized as follows: In \rsect{S1557} we recall some results on the Schur complement $\Hu{j}/\Hu{j-1}$ called \tSp{} $(\Spu{j})_{j=0}^\kappa$ of a sequence $(s_j)_{j=0}^\kappa$. In \rsect{S0820} the \tDSp{} given by the sequence $[\seq{\DSpl{k}}{k}{0}{\infi},\seq{\DSpm{k}}{k}{0}{\infi}]$ is recalled. Interrelations between this pair of sequences and the \tSp{} and vice versa are discussed. In \rsect{S1538}, via a set of nonnegative column pairs $\tmat{\phi(z)\\\psi(z)}$, the parametrization of the solution set in the non-degenerate case is rewritten in terms of a linear fractional transformation. Moreover, the representation of the resolvent matrix corresponding to a matricial truncated Stieltjes moment problem with matrix polynomials orthogonal on $\ra$ and matrix polynomials of the second  kind is recalled. In \rsect{S1649} we consider the notion and main results concerning the \tKt{}. \rsect{S1011} is devoted to the \tDSp{} $[\seq{\DSpol{k}{\ell}}{k}{0}{\infi},\seq{\DSpom{k}{\ell}}{k}{0}{\infi}]$ of the \taKta{\ell}{  $\seqsinf $}. A representation of the resolvent matrix of the matricial truncated Stieltjes moment problem in terms of the \tKt{ed} moment sequence is obtained in \rsect{S0834}. In \rsect{S1945}, interrelations between the polynomials $P^{(1)}_{k,j}$ and $Q^{(1)}_{k,j}$ constructed from the \tKt{ed} moment sequence with the polynomials $P_{k,j}$ and $Q_{k,j}$ corresponding to the original moment sequence are presented. A \tnnH{} measure, for which the matrix polynomials of the second kind are orthogonal matrix polynomials, is attained.

\section{\tSp{}}\label{S1557}
 With later applications to the matrix version of the \Stieltjes{} moment problem in mind, a particular inner parametrization, called \tSp{}, for matrix sequences was developed in~\zita{MR3014201}. First we are going to recall the definition of the \tSp{} of a sequence of complex \tpqa{matrices}. To prepare this notion, we need some further matrices built from the given data. If \(A\in\Cpq\), then a unique matrix $G\in\Cqp$ exists which satisfies the four equations $AGA=A$, $GAG=G$, $(AG)^\ad=AG$ and $(GA)^\ad=GA$. This matrix $G$ is called the \noti{Moore-Penrose inverse of $A$}{Moore-Penrose inverse} and is denoted by \(A^\MP\)\index{\(A^\MP\)}

 Let \(\kappa\in\NOinf \) and let \(\seqska \) be a sequence of complex \tpqa{matrices}. Then, let
\begin{align}\label{yz}
 \yuu{\ell}{m}&\defg
 \bMat
  \su{\ell}\\
  \su{\ell+1}\\
  \vdots\\
  \su{m}
 \eMat&
 &\text{and}&
 \zuu{\ell}{m}&\defg\brow\su{\ell},\su{\ell+1},\dotsc,\su{m}\erow
\end{align}
 \index{y@\(\yuu{\ell}{m}\)}\index{z@\(\zuu{\ell}{m}\)}for all \(\ell,m\in\NO\) with \(\ell\leq m\leq\kappa\). We use the notation
\begin{align}\label{L}
 \Lu{0}&\defg\su{0}&
 &\text{and}&
 \Lu{n}&\defg\su{2n}-\zuu{n}{2n-1}\Hu{n-1}^\MP\yuu{n}{2n-1}
\end{align}
 \index{l@\(\Lu{n}\)}for all \(n\in\N\) with \(2n\leq\kappa\), and the notation
\begin{align}\label{La}
 \Lau{0}&\defg\su{1}&
 &\text{and}&
 \Lau{n}&\defg\su{2n+1}-\zuu{n+1}{2n}\Ku{n-1}^\MP\yuu{n+1}{2n}
\end{align}
 \index{l@\(\Lau{n}\)}for all \(n\in\N\) with \(2n+1\leq\kappa\). Observe that for \(n\geq1\) the matrix \(\Lu{n}\) is the Schur complement \(\Hu{n}/\Hu{n-1}\) of \(\Hu{n-1}\) in the block \Hankel{} matrix \(\Hu{n}=\tmat{\Hu{n-1}&\yuu{n}{2n-1}\\\zuu{n}{2n-1}&\su{2n}}\) corresponding to the sequence \(\seqska \), whereas the matrix \(\Lau{n}\) is the Schur complement \(\Ku{n}/\Ku{n-1}\) of \(\Ku{n-1}\) in the block \Hankel{} matrix \(\Ku{n}=\tmat{\Ku{n-1}&\yuu{n+1}{2n}\\\zuu{n+1}{2n}&\su{2n+1}}\)
 corresponding to the sequence \(\seq{\sau{j}}{j}{0}{\kappa-1}\) defined by
\bgl{s+}
 \sau{j}
 \defg\su{j+1}.
\eg
 \index{\(\sau{j}\)}The sequence \(\seq{\sau{j}}{j}{0}{\kappa-1}\) coincides with the right-sided version of the sequence in~\zitaa{MR3014201}{\cdefnp{2.1}{217}} for \(\alpha=0\).

 Now we are able to recall the notion of \tSp{} which was introduced in~\zitaa{MR3014201}{\cdefnp{4.2}{223}} as \emph{right-sided \(\alpha\)\nobreakdash-\Stieltjes{} parametrization} denoted by \(\seq{Q_j}{j}{0}{\kappa}\) in a more general context related to a semi-infinite interval \([\alpha,\infp)\). There one can find further details.

\bdefil{D1021}
 Let \(\kappa\in\NOinf \) and let \(\seqska \) be a sequence of complex \tpqa{matrices}. Then the sequence \(\seq{\Spu{j}}{j}{0}{\kappa}\)\index{q@\(\seq{\Spu{j}}{j}{0}{\kappa}\)} given by \(\Spu{2k}\defg\Lu{k}\) for all \(k\in\NO\) with \(2k\leq\kappa\), and by \(\Spu{2k+1}\defg\Lau{k}\) for all \(k\in\NO\) with \(2k+1\leq\kappa\) is called the \noti{\tSpa{\(\seqska \)}}{Stieltjes!parametrization}.
\end{defn}

 There is a one-to-one correspondence between a sequence \(\seqska \) of complex \tpqa{matrices} and its
 \tSp{} \(\seq{\Spu{j}}{j}{0}{\kappa}\) (see~\zitaa{MR3014201}{\cremp{4.3}{224}}). In particular, the original sequence can be explicitly reconstructed from its \tSp{}. Let \symba{\nul{A}}{n} be the null space of a complex matrix \(A\).

\bpropnl{see~\zitaa{MR3014201}{\cthmp{4.12(b)}{225}}}{P0926}
 Let \(\kappa\in\NOinf \) and let \(\seqska \) be a sequence of complex \tqqa{matrices} with \tSp{} \((\Spu{j})_{j=0}^\kappa\). Then \(\seqska \in\Kggqu{\kappa}\) if and only if \(\Spu{j}\in\Cggq\) for all \(j\in\mn{0}{\kappa}\) and \(\nul{\Spu{j}}\subseteq\nul{\Spu{j+1}}\) for all \(j\in\mn{0}{\kappa-2}\).
\eprop

\bpropnl{see~\zitaa{MR3014201}{\cthmp{4.12(c)}{225}}}{P0937}
 Let \(\kappa\in\NOinf \) and let \(\seqska \) be a sequence of complex \tqqa{matrices} with \tSp{} \((\Spu{j})_{j=0}^\kappa\). Then \(\seqska \in\Kggequ{\kappa}\) if and only if \(\Spu{j}\in\Cggq\) for all \(j\in\mn{0}{\kappa}\) and \(\nul{\Spu{j}}\subseteq\nul{\Spu{j+1}}\) for all \(j\in\mn{0}{\kappa-1}\).
\eprop

 Now we introduce an important subclass of the class of \tSnnd{} sequences. More precisely, we turn our attention to some subclass of \(\Kggqu{\kappa}\), which is characterized by stronger positivity properties. For all \(n\in\NO\), let \(\Kgqu{2n+1}\)\index{k@\(\Kgqu{2n+1}\)} be the set of all sequences \(\seq{\su{j}}{j}{0}{2n+1}\) of complex \tqqa{matrices}, such that \(\Hu{n}\) and \(\Ku{n}\) are both \tpH{}. For all \(n\in\NO\), let \(\Kgqu{2n}\)\index{k@\(\Kgqu{2n}\)} be the set of all sequences \(\seq{\su{j}}{j}{0}{2n}\) of complex \tqqa{matrices}, such that \(\Hu{n}\) is \tpH{} and, in the case \(n\geq1\), furthermore \(\Ku{n-1}\) is \tpH{}. Let \(\Kgqu{\infty}\)\index{k@\(\Kgqinf\)} be the set of all sequences \(\seqsinf \) of complex \tqqa{matrices} such that \(\seq{\su{j}}{j}{0}{m}\) belongs to \(\Kgqu{m}\) for all \(m\in\NO\). For all \(\kappa\in\NOinf \), we call a sequence \(\seqska \) \noti{\tSpd}{Stieltjes!positive definite} if it belongs to \(\Kgqu{\kappa}\). For \(\kappa\in\NOinf \), we have \(\Kgqkappa\subseteq\Kggeqkappa\subseteq\Kggqkappa\) (see~\zitaa{114arxiv}{\cpropp{3.8}{12}}). In view of \rtheoss{T1121}{T1122} we obtain then:

\begin{rem}\label{R1524}
 Let \(\kappa\in\NOinf \) and let \(\seqska \in\Kgqu{\kappa}\). Then, \(\Mggqaag{\ra}{\seqska }\neq\emptyset\) and, in the case \(\kappa<\infi\), furthermore \(\Mggqaakg{\ra}{\seqska }\neq\emptyset\).
\end{rem}

 If \(m\in\NO\) and \(\seqs{m}\in\Kgqu{m}\), the associated truncated moment problems~\mproblem{\ra}{m}{=} and~\mproblem{\ra}{m}{\leq} have infinitely many solutions. Since every principal submatrix of a \tpH{} matrix is again \tpH{}, we can easily see:
\bremal{R2347}
 Let \(\kappa\in\Ninf\) and let \(\seqska\in\Kgqkappa\), then \(\seq{\sau{j}}{j}{0}{\kappa-1}\in\Kgqu{\kappa-1}\).
\erema

\bpropnl{see~\zitaa{MR3014201}{\cthmp{4.12(d)}{225}}}{T1337}
 Let \(\kappa\in\NOinf \) and let \(\seqska \) be a sequence of complex \tqqa{matrices} with \tSp{} \((\Spu{j})_{j=0}^\kappa\). Then \(\seqska \in\Kgqu{\kappa}\) if and only if \(\Spu{j}\in\Cgq\) for all \(j\in\mn{0}{\kappa}\).
\eprop

 Based on the matrices defined via \eqref{L} and \eqref{La}, we now introduce a further important subclass of $\Kggqkappa$. Let $m\in\NO$ and let $\seqs{m}\in\Kggqu{m}$. Then $\seqs{m}$ is called \notii{completely degenerate} if $\Lu{n}=\Oqq$ in the case \(m=2n\) with some \(n\in\NO\) or if $\Lau{n}=\Oqq$ in the case \(m=2n+1\) with some \(n\in\NO\). The set \symba{\Kggdqu{m}}{k} of all completely degenerate sequences belonging to $\Kggqu{m}$ is a subset of $\Kggequ{m}$ (see~\zitaa{MR3014201}{\cpropp{5.9}{231}}). The moment problem~\mproblem{\ra}{m}{=} has a unique solution if and only if \(\seqs{m}\) belongs to \(\Kggdqu{m}\) (see~\zitaa{142arxiv}{\cthmp{13.3}{53}}).

\bpropnl{cf.~\zitaa{MR3014201}{\cpropp{5.3}{229}}}{111.P5-3}
 Let $m\in\NO$ and $\seqs{m}\in\Kggqu{m}$ with \tSp{} \((\Spu{j})_{j=0}^m\). Then $\seqs{m}\in\Kggdqu{m}$ if and only if $\Spu{m}=\Oqq$.
\eprop

 If $\seqs{m}\in\Kggequ{m}$ with \tSp{} \((\Spu{j})_{j=0}^m\), then from~\zitaa{MR2735313}{\clemmss{4.15}{4.16}} one can easily see that $\seqs{m}$ belongs to $\Kggdqu{m}$ if and only if there is some $\ell\in\mn{0}{m}$ such that $\Spu{\ell}=\Oqq$.

 Let $\seqsinf\in\Kggqinf$. Then \(\seqsinf\) is said to be \notii{completely degenerate} if there is some $m\in\NO$ such that $\seqs{m}$ is a \tcdSnnd{} sequence. By \symba{\Kggdqinf}{k} we denote the set of all \tcdSnnd{} sequences $\seqsinf $ of complex \tqqa{matrices}.

\bpropnl{cf.~\zitaa{MR3014201}{\ccorp{5.4}{230}}}{111.C5-4}
 Let $\seqsinf\in\Kggqinf$ with \tSp{} \((\Spu{j})_{j=0}^\infi\). Then $\seqsinf\in\Kggdqinf$ if and only if there exists some \(m\in\NO\) with $\Spu{m}=\Oqq$.
\eprop

 The sequence $\seqsinf $  is called \notii{completely degenerate of order $m$} if $\seqs{m}$ is completely degenerate. By \symba{\Kggdoq{m}}{k} we denote the set of all \tSnnd{} sequences $\seqsinf $ from $\Cqq$ which are completely degenerate of order $m$. If $m\in\NO$ and $\seqsinf \in\Kggdoq{m}$, then $\seqs{\ell}\in\Kggdqu{\ell}$ for each $\ell\in\minf{m}$ (cf.~\zitaa{MR3014201}{\clemp{5.5}{230}}).

\bdefil{D0909}
 Let \(m\in\NO\) and let \(\seqs{m}\) be a sequence of complex \tpqa{matrices}. Let the sequence \(\seq{\Spl{j}}{j}{0}{\infi}\) be given by
 \[
  \Spl{j}
  \defg
  \begin{cases}
   \Spu{j}\ifa{j\leq m}\\
   \Opq\ifa{j>m}
  \end{cases}.
 \]
 Then, we call the unique sequence \(\seqzinf\)\index{$\seqzinf$} with \tSp{} \(\seq{\Spl{j}}{j}{0}{\infi}\) the \notii{\tzexto{\(\seqs{m}\)}}.
\edefi

\blemml{L0916}
 Let \(m\in\NO\) and let \(\seqs{m}\in\Kggequ{m}\). Denote by \(\seqzinf\) the \tzexto{\(\seqs{m}\)}. Then, \(\zext{s}{j}=\su{j}\) for all \(j\in\mn{0}{m}\) and \(\seqzinf\in\Kggdoq{m+1}\).
\elemm
\bproof
 Denote by \(\seq{\Spu{j}}{j}{0}{m}\) the \tSpa{$\seqs{m}$} and by \(\seq{\Spl{j}}{j}{0}{\infi}\) the \tSpa{$\seqzinf$}. Since \(\Spl{j}=\Spu{j}\) holds true for all \(j\in\mn{0}{m}\), we have \(\zext{s}{j}=\su{j}\) for all \(j\in\mn{0}{m}\). Using \rpropss{P0937}{P0926} we can conclude furthermore \(\seqzinf\in\Kggqinf\). In view of \(\Spl{m+1}=\Oqq\) and \rprop{111.P5-3}, thus \(\seqzinf\in\Kggdoq{m+1}\) follows.
\eproof

\blemml{L0912}
 Let \(m\in\NO\) and let \(\seqsinf\in\Kggdoq{m+1}\). Denote by \(\seqzinf\) the \tzexto{\(\seqs{m}\)}. Then, \(\zext{s}{j}=\su{j}\) for all \(j\in\NO\).
\elemm
\bproof
 Denote by \(\seq{\Spu{j}}{j}{0}{\infi}\) the \tSpa{$\seqsinf$} and by \(\seq{\Spl{j}}{j}{0}{\infi}\) the \tSpa{$\seqzinf$}. Since \(\seq{\Spu{j}}{j}{0}{m}\) is then the \tSpa{$\seqs{m}$}, we have by definition \(\Spl{j}=\Spu{j}\) for all \(j\in\mn{0}{m}\). Because of \(\seqsinf\in\Kggdoq{m+1}\), the sequence \(\seqs{m+1}\) belongs to \(\Kggdqu{m+1}\). Thus, \rprop{111.P5-3} yields \(\Spu{m+1}=\Oqq\). From \rprop{P0926}, we obtain then \(\Spu{j}=\Oqq\) for all \(j\in\mn{m+2}{\infi}\). By definition, we have furthermore \(\Spl{j}=\Oqq\) for all \(j\in\mn{m+1}{\infi}\). Hence, \(\Spl{j}=\Spu{j}\) for all \(j\in\mn{m+1}{\infi}\). We have shown that the \tSp{s} of $\seqsinf$ and $\seqzinf$ coincide, which completes the proof.
\eproof

\section{\tDSp{}}\label{S0820}
 In~\zita{MR2053150} Yu.~M.~\tDyukarev{} studied the moment problem \mproblem{\ra}{\infty}{=}. One of his main results (see~\zitaa{MR2053150}{\cthmp{8}{78}}) is a generalization of a classical criterion due to \Stieltjes{}~\zitas{MR1508159,MR1508160} for the indeterminacy of this moment problem. In order to find an appropriate matricial version of \Stieltjes{'} indeterminacy criterion Yu.~M.~\tDyukarev{} had to look for a convenient matricial generalization of the parameter sequences which \Stieltjes{} obtained from the consideration of particular continued fractions associated with the sequence \(\seqsinf \). In this way, Yu.~M.~\tDyukarev{} found an interesting inner parametrization of sequences belonging to \(\Kgqinf\). The main theme of this section is to recall some interrelations obtained in~\zitaa{MR3133464}{\S8} between Yu.~M.~\tDyukarev{'s} parametrization and the \tSp{} introduced in \rdefi{D1021}.

 The notations \(\Iq\)\index{i@\(\Iq\)} and \(\Opq\)\index{0@\(\Opq\)} stand for the identity matrix in \(\Cqq\) and for the zero matrix in \(\Cpq\), resp. If \(\kappa\in\NOinf \) and \(\seqska \in\Kgqkappa\), then the matrix \(\Hu{k}\) is \tpH{} and, in particular, invertible for all \(k\in\NO\) with \(2k\leq\kappa\), and the matrix \(\Ku{k}\) is \tpH{} and, in particular, invertible for all \(k\in\NO\) with \(2k+1\leq\kappa\). Let
\begin{align}\label{v}
 \vqu{0}&\defg\Iq&
 &\text{and}&
 \vqu{k}&\defg\bMat\Iq\\\Ouu{kq}{q}\eMat
\end{align}
\index{v@\(\vqu{k}\)}for all \(k\in\N\). The following construction of a pair of sequences of \tqqa{matrices} associated with a \tSpd{} sequence goes back to Yu.~M.~\tDyukarev{}~\zitaa{MR2053150}{p.~77}:

Let \(\kappa\in\NOinf \) and let \(\seqska \in\Kgqkappa\). Then let
\beql{M0}
 \DSpm{0}
 \defg\su{0}^\inv
\eeq
and, in the case \(\kappa\geq1\), let
\beql{L0}
 \DSpl{0}
 \defg\su{0}\su{1}^\inv\su{0}.
\eeq
Furthermore, let
\beql{Mk}
 \DSpm{k}
 \defg\vqu{k}^\ad\Hu{k}^\inv\vqu{k}-\vqu{k-1}^\ad\Hu{k-1}^\inv\vqu{k-1}
\eeq
\index{m@\(\DSpm{k}\)}for all \(k\in\N\) with \(2k\leq\kappa\), and let
\beql{Lk}
 \DSpl{k}
 \defg\yuu{0}{k}^\ad\Ku{k}^\inv\yuu{0}{k}-\yuu{0}{k-1}^\ad\Ku{k-1}^\inv\yuu{0}{k-1}
\eeq
\index{l@\(\DSpl{k}\)}for all \(k\in\N\) with \(2k+1\leq\kappa\).

Obviously, for all \(k\in\NO\) with \(2k\leq\kappa\), the matrix \(\DSpm{k}\)
only depends on the matrices \(\su{0},\dotsc,\su{2k}\), and, for all \(k\in\NO\) with \(2k+1\leq\kappa\), the matrix \(\DSpl{k}\) only depends on the matrices \(\su{0},\su{1},\dotsc,\su{2k+1}\).

\bdefil{D1455}
 Let \(\seqsinf \in\Kgqinf\), then the ordered pair \([\seq{\DSpl{k}}{k}{0}{\infty},\seq{\DSpm{k}}{k}{0}{\infty}]\)
 is called the \noti{\tDSp{}}{d@\tDSp{}} (shortly \noti{\tdsp{}}{d@\tDSp{}}) \emph{of \(\seqsinf \)}.
\edefi

 It should be mentioned that, for a given sequence \(\seqsinf \in\Kgqinf\), Yu.~M.~Dyukarev~\zita{MR2053150} treated the moment problem~\mproblem{\ra}{\infi}{=} by approximation through the sequence \((\mproblem{\ra}{k}{\leq})_{k\in\NO}\) of truncated moment problems. One of his central results~\zitaa{MR2053150}{\cthmp{7}{77}} shows that the resolvent matrices for the truncated moment problems can be multiplicatively decomposed into elementary factors which are determined by the corresponding first sections of the \tdspa{\(\seqsinf \)}.

 If \(\seqsinf \in\Kgqinf\) with \tdsp{} \([\seq{\DSpl{k}}{k}{0}{\infty},\seq{\DSpm{k}}{k}{0}{\infty}]\), then, in view of~\zitaa{MR2053150}{\cthmp{7}{77}}, the matrices \(\DSpl{k}\) and \(\DSpm{k}\) are \tpH{} for all \(k\in\NO\). According to~\zitaa{MR3133464}{\cpropp{8.26}{3923}}, every sequence \(\seqsinf \in\Kgqinf\) can be recursively reconstructed from its \tdsp{} \([\seq{\DSpl{k}}{k}{0}{\infty},\seq{\DSpm{k}}{k}{0}{\infty}]\). A similar result for the truncated \tdsp{} \([\seq{\DSpl{k}}{k}{0}{m-1},\seq{\DSpm{k}}{k}{0}{m}]\) (resp.\ \([\seq{\DSpl{k}}{k}{0}{m},\seq{\DSpm{k}}{k}{0}{m}]\)) was obtained in~\zitaa{MR3327132}{\cpropp{4.9}{68}}. Furthermore,~\zitaa{MR3133464}{\cpropp{8.27}{3924}} shows that each pair \([\seq{\DSpl{k}}{k}{0}{\infty},\seq{\DSpm{k}}{k}{0}{\infty}]\) of sequences of \tpH{} complex \tqqa{matrices} is the \tdsp{} of some sequence \(\seqsinf \in\Kgqinf\). Hence, the \tdsp{} establishes a one-to-one correspondence between \tSpd{} sequences \(\seqsinf \) and ordered pairs \([\seq{\DSpl{k}}{k}{0}{\infty},\seq{\DSpm{k}}{k}{0}{\infty}]\) of sequences of \tpH{} complex \tqqa{matrices}.

 In~\zitaa{MR3133464}{\cpropp{8.30}{3925}} it was shown that in the scalar case the \tdsp{} of a \tSpd{} sequence coincides with the classical parameters used by \Stieltjes{}~\zitas{MR1508159,MR1508160} to formulate his indeterminacy criterion. We mention that M.~G.~Kre{\u\i}n was able to find a mechanical interpretation for \Stieltjes{'} investigations on continued fractions (see Gantmacher/Kre{\u\i}n~\zitaa{MR0114338}{\canha{2}} or Akhiezer~\zitaa{MR0184042}{Appendix}). Against to the background of his mechanical interpretation M.~G.~Kre{\u\i}n divided \Stieltjes{'} original parameters into two groups which play the roles of lengths and masses, respectively.
 Now we want to recall the concrete definition of these parameters (see Kre{\u\i}n/Nudel{\('\)}man~\zitaa{MR0458081}{\cchap{V}, \cform{(6.1)}}) and their connection to the \tdsp{}.

\bdefil{D0834}
 Let \(\seqsinf \in\Kguu{1}{\infty}\) and let \(\detHu{n}\defg\det\Hu{n}\)\index{d@\(\detHu{n}\)} and \(\detKu{n}\defg\det\Ku{n}\)\index{d@\(\detKu{n}\)} for all \(n\in\NO\). Let
 \begin{align*}
  \spl{k}&\defg\frac{\detHu{k}^2}{\detKu{k}\detKu{k-1}}&
  &\text{and}&
  \spm{k}&\defg\frac{(\detKu{k-1})^2}{\detHu{k}\detHu{k-1}}
 \end{align*}
 \index{l@\(\spl{k}\)}\index{m@\(\spm{k}\)}for all \(k\in\NO\), where \(\detHu{-1}\defg1\) and \(\detKu{-1}\defg1\). Then the ordered pair \([\seq{\spl{k}}{k}{0}{\infty},\seq{\spm{k}}{k}{0}{\infty}]\) is called the \noti{\tKSp{} of \(\seqsinf \)}{k@\tKSp{}}.
\edefi

\bpropnl{see~\zitaa{MR3133464}{\cpropp{8.30}{3925}}}{P0857}
 Let \(\seqsinf \in\Kguu{1}{\infty}\) with \tKSp{} \([\seq{\spl{k}}{k}{0}{\infty},\seq{\spm{k}}{k}{0}{\infty}]\) and \tdsp{} \([\seq{\DSpl{k}}{k}{0}{\infty},\seq{\DSpm{k}}{k}{0}{\infty}]\). Then \(\spl{k}=\DSpl{k}\) and \(\spm{k}=\DSpm{k}\) for all \(k\in\NO\).
\eprop

Now we recall the connection between the \tdsp{} of a \tSpd{} sequence and its \tSp{}. If \(\seqsinf \in\Kgqinf\) with \tSp{} \(\seq{\Spu{j}}{j}{0}{\infty}\), then, in view of \rprop{T1337}, the matrices \(\Spu{j}\) are \tpH{} and, in particular, invertible for all \(j\in\NO\).

\bthmnl{see~\zitaa{MR3133464}{\cthmp{8.22}{3921}}}{T1513}
 Let \(\seqsinf \in\Kgqinf\) with \tSp{} \(\seq{\Spu{j}}{j}{0}{\infty}\) and \tdsp{} \([\seq{\DSpl{k}}{k}{0}{\infty},\seq{\DSpm{k}}{k}{0}{\infty}]\). Then
 \[
  \DSpl{k}
  =\rk*{\rprod_{j=0}^k\Spu{2j}\Spu{2j+1}^\inv}\Spu{2k+1}\rk*{\rprod_{j=0}^k\Spu{2j}\Spu{2j+1}^\inv}^\ad
 \]
 and
 \[
  \DSpm{k}
  =
  \begin{cases}
    \Spu{0}^\inv\incase{k=0}\\
    \rk*{\rprod_{j=0}^{k-1}\Spu{2j}^\inv\Spu{2j+1}}\Spu{2k}^\inv\rk*{\rprod_{j=0}^{k-1}\Spu{2j}^\inv\Spu{2j+1}}^\ad\incase{k\geq1}
  \end{cases}
 \]
 for all \(k\in\NO\).
\etheo
\bthmnl{see~\zitaa{MR3133464}{\cthmp{8.24}{3923}} and
\zitaa{MR3327132}{Corollary 4.10, p. 72}}{T1523}
 Let \(\seqsinf \in\Kgqinf\) with \tSp{} \(\seq{\Spu{j}}{j}{0}{\infty}\) and \tdsp{} \([\seq{\DSpl{k}}{k}{0}{\infty},\seq{\DSpm{k}}{k}{0}{\infty}]\). Then
 \[
  \Spu{2k}
  =
  \begin{cases}
    \DSpm{0}^\inv\incase{k=0}\\
    \rk*{\rprod_{j=0}^{k-1}\DSpm{j}\DSpl{j}}^\invad\DSpm{k}^\inv\rk*{\rprod_{j=0}^{k-1}\DSpm{j}\DSpl{j}}^\inv\incase{k\geq1}
  \end{cases}
 \]
 and
 \[
  \Spu{2k+1}
  =\rk*{\rprod_{j=0}^k\DSpm{j}\DSpl{j}}^\invad\DSpl{k}\rk*{\rprod_{j=0}^k\DSpm{j}\DSpl{j}}^\inv
 \]
 for all \(k\in\NO\).
\etheo

%

\section{Parametrization of all solutions in the non-degenerate case}\label{S1538}
 In this section we recall a parametrization of the solution set of the moment problem \mproblem{\ra}{m}{\leq} for the so-called non-degenerate case. Let \(\SFq\)\index{s@\(\SFq\)} be the set of all holomorphic functions \(F\colon\Cs\to\Cqq\) which satisfy the following conditions:
\bAeqi{0}
 \item The matrix \(\im F(w)\) is \tnnH{} for all \(w\in\C\) with \(\im w>0\).
 \item The matrix \(F(x)\) is \tnnH{} for all \(x\in(-\infty,0)\).
\eAeqi
 Further, let \(\SOFq\)\index{s@\(\SOFq\)} be the set of all \(S\in\SFq\) such that
\[
 \sup_{y\in[1,\infp)}y\normAE{S(\ii y)}
 <\infp,
\]
 where \(\normaE{\cdot}\)\index{e@\(\normaE{\cdot}\)} is the Euclidean matrix norm. We have the following integral representation for functions belonging to \(\SOFq\):
\bthmnl{cf.~\zitaa{142arxiv}{\cthmp{5.1}{19}}}{T1434}
 Let \(S\colon\Cs\to \Cqq\).
 \benui
  \il{T1434.a} If \(S\in \SOFq\), then there exists a unique \tnnH{} measure \(\sigma\in\Mggqa{\ra}\) such that
  \begin{equation}\label{T1434.1}
   S(z)
   =\int_{\ra} \frac{1}{t-z} \sigma(\dif t)
  \end{equation}
  for all \(z\in\Cs\).
  \il{T1434.b} If there exists a \tnnH{} measure \(\sigma\in\Mggqa{\ra}\) such that \(S\) can be represented via \eqref{T1434.1} for all \(z\in\Cs\), then \(S\) belongs to \(\SOFq\).
 \eenui
\etheo
 If \(\sigma\) is a measure belonging to \(\Mggqa{\ra}\), then we will call the matrix-valued function \(S\colon\Cs\to\Cqq\) which is given for all \(z\in\Cs\) by \eqref{T1434.1} the \noti{\tSta{\(\sigma\)}}{Stieltjes!transform} and write \(\sttrhl{\sigma}\) for \(S\). If \(S\in\SOFq\), then the unique measure \(\sigma\) which belongs to \(\Mggqa{\ra}\) and which fulfills \eqref{T1434.1} for all \(z\in\Cs\) is said to be the \noti{\tSma{\(S\)}}{Stieltjes!measure}.

In view of \rtheo{T1434}, the moment problem~\mproblem{\ra}{m}{\leq} can be reformulated:
\begin{quote}
 Let \(m\in\NO\) and let \(\seq{\su{j}}{j}{0}{m}\) be a sequence of complex \tqqa{matrices}. Describe the set \(\SOFqakg{\seq{\su{j}}{j}{0}{m}}\)\index{s@\(\SOFqakg{\seq{\su{j}}{j}{0}{m}}\)} of all \(S\in\SOFq\) with \tSm{} belonging to \(\Mggqaakg{\ra}{\seq{\su{j}}{j}{0}{m}}\).
\end{quote}

Let the \taaa{2q}{2q}{signature} matrices \(\Jreq\) and \(\Jimq\)
 associated with the real and imaginary part be defined by
\begin{align*}
 \Jreq
 &\defg
 \bMat
  \Oqq&\Iq\\
  \Iq&\Oqq
 \eMat&
 &\text{and}&
 \Jimq
 &\defg
 \bMat
  \Oqq&\ii\Iq\\
  -\ii\Iq&\Oqq
 \eMat.
\end{align*}
\index{j@\(\Jreq\)}\index{j@\(\Jimq\)}Let \(\SPq\)\index{s@\(\SPq\)} be the set of all ordered pairs \((\phi,\psi)\) of \(\Cqq\)\nobreakdash-valued functions \(\phi\) and \(\psi\) which are meromorphic in \(\Cs\) and for which there exists a discrete subset \(\Dpq\) of \(\Cs\) such that the following conditions are fulfilled:
\bAeqi{2}
 \item The functions \(\phi\) and \(\psi\) are both holomorphic in \(\C\setminus(\ra\cup\Dpq)\).
 \item \(\rank\tmat{\phi(z)\\\psi(z)}=q\) for all \(z\in\C\setminus(\ra\cup\Dpq)\).
 \item For all \(z\in\C\setminus(\R\cup\Dpq)\),
 \[
  \frac{1}{\im z}
  \bMat
   \phi(z)\\
   \psi(z)
  \eMat^\ad\Jimq
  \bMat
   \phi(z)\\
   \psi(z)
  \eMat
  \in\Cggq.
 \]
 \item For all \(z\in\C\setminus\Dpq\) with \(\re z<0\),
 \[
  \bMat
   \phi(z)\\
   \psi(z)
  \eMat^\ad\Jreq
  \bMat
   \phi(z)\\
   \psi(z)
  \eMat
  \in\Cggq.
 \]
\eAeqi
Two pairs \((\phi_1,\psi_1),(\phi_2,\psi_2)\in\SPq\) are called \noti{equivalent}{pairs!equivalent}, if there exists a \(\Cqq\)\nobreakdash-valued function \(\eta\) which is meromorphic in \(\Cs\) such that \(\det\eta\) does not vanishing identically and
\[
 \bMat
  \phi_2\\
  \psi_2
 \eMat
 =
 \bMat
  \phi_1\\
  \psi_1
 \eMat\eta.
\]
For all \(n\in\NO\), let
\[
 \Tqu{n}
 \defg\matauuo{\Kronu{j,k+1}\Iq}{j,k}{0}{n}
\]
\index{t@\(\Tqu{n}\)}and let \(\Rqu{n}\colon\C\to\Coo{(n+1)q}{(n+1)q}\)\index{r@\(\Rqu{n}\)} be defined by
\[
 \Rqua{n}{z}
 \defg(\Iu{(n+1)q}-z\Tqu{n})^\inv.
\]
Let \(\kappa\in\NOinf \) and let \(\seqska \) be a sequence of complex \tqqa{matrices}. Then let
\begin{align*}
 \uu{0}&\defg\Oqq&
 &\text{and}&
 \uu{k}&\defg\bMat\Oqq\\-\yuu{0}{k-1}\eMat
\end{align*}
\index{u@\(\uu{k}\)}for all \(k\in\mn{1}{\kappa+1}\).
 Now we suppose that \(\seqska \)
  belongs to \(\Kgqkappa\). Let
  \(\pau{n}\colon\C\to\Cqq\)\index{a@\(\pau{n}\)} and
  \(\pcu{n}\colon\C\to\Cqq\)\index{c@\(\pcu{n}\)} be defined by
\begin{align}
 \paua{n}{z}&\defg\Iq-z\uu{n}^\ad\ek*{\Rqua{n}{\ko z}}^\ad
 \Hu{n}^\inv\vqu{n} \label{alp0}
\shortintertext{and}
 \pcua{n}{z}&\defg-z\vqu{n}^\ad\ek*{\Rqua{n}{\ko z}}^\ad
 \Hu{n}^\inv\vqu{n}\label{bet0}
\end{align}
for all \(n\in\NO\) with \(2n\leq\kappa\), and let \(\pbu{n}\colon\C\to\Cqq\)\index{b@\(\pbu{n}\)} and \(\pdu{n}\colon\C\to\Cqq\)\index{d@\(\pdu{n}\)} be defined by
\begin{align}
 \pbua{n}{z}&\defg
 \begin{cases}
  \Oqq\incase{n=0}\\
  \yuu{0}{n-1}^\ad\ek*{\Rqua{n-1}{\ko z}}^\ad\Ku{n-1}^\inv
  \yuu{0}{n-1}\incase{n\geq1} \label{gam0}
 \end{cases}
\shortintertext{and}
 \pdua{n}{z}&\defg
 \begin{cases}
  \Iq\incase{n=0}\\
  \Iq+z\vqu{n-1}^\ad\ek*{\Rqua{n-1}{\ko z}}^\ad\Ku{n-1}^\inv
  \yuu{0}{n-1}\incase{n\geq1} \label{del0}
 \end{cases}
\end{align}
 for all \(n\in\NO\) with \(2n-1\leq\kappa\). Then
\begin{align*}
 \paua{0}{z}&=\Iq,&&&\pbua{0}{z}&=\Oqq\\
 \pcua{0}{z}&=-z\su{0}^\inv,&&\text{and}&\pdua{0}{z}&=\Iq
\end{align*}
for all \(z\in\C\). Let
\begin{equation} \label{Uu2n}
 \Uu{2n}
 \defg
 \bMat
  \pau{n}&\pbu{n}\\
  \pcu{n}&\pdu{n}\\
 \eMat
\end{equation}
\index{u@\(\Uu{2n}\)}for all \(n\in\NO\) with \(2n\leq\kappa\), and let
\begin{equation} \label{Uu2n1}
 \Uu{2n+1}
 \defg
 \bMat
  \pau{n}&\pbu{n+1}\\
  \pcu{n}&\pdu{n+1}\\
 \eMat
\end{equation}
 \index{u@\(\Uu{2n+1}\)}for all \(n\in\NO\) with
\(2n+1\leq\kappa\).

\bdefil{D1151}
 Let \(\seqsinf\in\Kgqinf\). Then \symba{\seq{\Uu{m}}{m}{0}{\infi}}{u} is called the \noti{\tsqDmpo{\(\seqsinf\)}}{d@\tqDyump{}}.
\edefi

 Now we are able to recall a parametrization of
the set \(\SOFqakg{\seq{\su{j}}{j}{0}{m}}\) with pairs belonging to
\(\SPq\):
 \bthmnl{see~\zitaa{MR3324594}{\cthmp{3.2}{9}} or~\zitas{Dyu81,MR1699439}}{T1327}
 Let \(m\in\NO\) and let \(\seq{\su{j}}{j}{0}{m}\in\Kgqu{m}\). Let \(\Uu{m}=\tmat{A_m&B_m\\C_m&D_m}\) be the \tqqa{block} representation of \(\Uu{m}\). Then:
 \benui
  \il{T1327.a} Let \((\phi,\psi)\in\SPq\). Then \(\det(C_m\phi+D_m\psi)\) does not vanish identically in \(\Cs\) and
  \[
   (A_m\phi+B_m\psi)(C_m\phi+D_m\psi)^\inv
   \in\SOFqakg{\seq{\su{j}}{j}{0}{m}}.
  \]
  \item Let \((\phi_1,\psi_1),(\phi_2,\psi_2)\in\SPq\) be such that
  \[
   (A_m\phi_1+B_m\psi_1)(C_m\phi_1+D_m\psi_1)^\inv
   =(A_m\phi_2+B_m\psi_2)(C_m\phi_2+D_m\psi_2)^\inv.
  \]
  Then, \((\phi_1,\psi_1)\) and \((\phi_2,\psi_2)\) are equivalent.
  \item Let \(S\in\SOFqakg{\seq{\su{j}}{j}{0}{m}}\). Then, there exists a pair \((\phi,\psi)\in\SPq\) such that
  \[
   S
   =(A_m\phi+B_m\psi)(C_m\phi+D_m\psi)^\inv.
  \]
 \eenui
\etheo

If \(\kappa\in\NOinf \) and \(\seqska \in\Kgqkappa\), then let \(\Mmu{k}\colon\C\to\Coo{2q}{2q}\)\index{m@\(\Mmu{k}\)} be defined by
\begin{equation}
 \Mmua{k}{z}
 \defg
 \bMat
  \Iq&\Oqq\\
  -z\DSpm{k}&\Iq
  \eMat
 \label{ma11}
\end{equation}
for all \(k\in\NO\) with \(2k\leq\kappa\), and let
\begin{equation}
 \Llu{k}
 \defg
 \bMat
  \Iq&\DSpl{k}\\
  \Oqq&\Iq
  \eMat
 \label{la11}
\end{equation}
\index{l@\(\Llu{k}\)}for all \(k\in\NO\) with \(2k+1\leq\kappa\).
 We have the following factorization of \(\Uu{m}\) in a product of
  complex \taaa{2q}{2q}{matrix} polynomials of degree \(1\):
\bpropnl{see~\zitaa{MR3324594}{\ceqqp{(4.11)}{(4.12)}{11}} or~\zitaa{MR2053150}{\cthmp{7}{77}}}{P1033}
 Let \(\kappa\in\NOinf \) and let
 \(\seqska \in\Kgqkappa\), then
 \begin{align*}
  \Uu{0}&=\Mmu{0}&
  &\text{and}&
  \Uu{2n}&=\rk*{\rprod_{k=0}^{n-1}\Mmu{k}\Llu{k}}\Mmu{n}
 \end{align*}
 for all \(n\in\N\) with \(2n\leq\kappa\), and
 \[
  \Uu{2n+1}
  =\rprod_{k=0}^n(\Mmu{k}\Llu{k})
 \]
 for all \(n\in\NO\) with \(2n+1\leq\kappa\).
\eprop


\bnotal{N1104}
 Let \(\kappa\in\NOinf \) and let
\(\seqska \in\Kgqkappa\). Let
\(\ophu{n}\colon\C\to\Cqq\)\index{p@\(\ophu{n}\)} and
\(\sphu{n}\colon\C\to\Cqq\)\index{q@\(\sphu{n}\)} be defined by
\begin{align*}
 \ophua{n}{z}&\defg
 \begin{cases}
  \Iq\incase{n=0}\\
  \vqu{n}^\ad\ek*{\Rqua{n}{\ko z}}^\ad\tmat{-\Hu{n-1}^\inv\yuu{n}{2n-1}\\\Iq}\incase{n\geq1}
 \end{cases}
\shortintertext{and}
 \sphua{n}{z}&\defg
 \begin{cases}
  \Oqq\incase{n=0}\\
  -\uu{n}^\ad\ek*{\Rqua{n}{\ko z}}^\ad\tmat{-\Hu{n-1}^\inv\yuu{n}{2n-1}\\\Iq}\incase{n\geq1}
 \end{cases}
\end{align*}
for all \(n\in\NO\) with \(2n-1\leq\kappa\), and let \(\opku{n}\colon\C\to\Cqq\)\index{p@\(\opku{n}\)} and \(\spku{n}\colon\C\to\Cqq\)\index{q@\(\spku{n}\)} be defined by
\begin{align*}
 \opkua{n}{z}&\defg
 \begin{cases}
  \Iq\incase{n=0}\\
  \vqu{n}^\ad\ek*{\Rqua{n}{\ko z}}^\ad\tmat{-\Ku{n-1}^\inv\yuu{n+1}{2n}\\\Iq}\incase{n\geq1}
 \end{cases}
\shortintertext{and}
 \spkua{n}{z}&\defg
 \begin{cases}
  \su{0}\incase{n=0}\\
  \zuu{0}{n}\ek*{\Rqua{n}{\ko z}}^\ad\tmat{-\Ku{n-1}^\inv\yuu{n+1}{2n}\\\Iq}\incase{n\geq1}
 \end{cases}
\end{align*}
for all \(n\in\NO\) with \(2n\leq\kappa\).
\enota

\bdefil{D1137}
 Let \(\seqsinf\in\Kgqinf\). Then \(\PQPQ\) is called the \noti{\tStiqosoompo{\(\seqsinf\)}}{s@\tStiqosoomp{}}.
\edefi

\bremal{R0001}
 For all \(n\in\NO\) with \(2n-1\leq\kappa\), the functions \(\ophu{n}\) and \(\sphu{n}\) are \tqqa{matrix} polynomials, where \(\ophu{n}\) has degree \(n\) and leading coefficient \(\Iq\). Similarly, for all \(n\in\NO\) with \(2n\leq\kappa\), the functions \(\opku{n}\) and \(\spku{n}\) are \tqqa{matrix} polynomials, where \(\opku{n}\) has degree \(n\) and leading coefficient \(\Iq\). Furthermore, if \(\sigmah\in\Mggqaag{\ra}{\seqska }\), we have the orthogonality relations
\[
\int_{\ra}\ek*{\ophua{m}{t}}^\ad\sigmah(\dif t)\ek*{\ophua{n}{t}}
 =
 \begin{cases}
   \Oqq\incase{m\neq n}\\
  \Lu{n}\incase{m=n}
 \end{cases}
\]
for all \(m,n\in\NO\) with \(2n-1\leq\kappa\) and \(2m-1\leq\kappa\), and
\[
 \int_{\ra}\ek*{\opkua{m}{t}}^\ad\sigmak(\dif t)\ek*{\opkua{n}{t}}
 =
 \begin{cases}
  \Oqq\incase{m\neq n}\\
  \Lau{n}\incase{m=n}
 \end{cases}
%
\]
for all \(m,n\in\NO\) with \(2n\leq\kappa\) and \(2m\leq\kappa\), where \(\sigmak\colon\Bra\to\Cggq\)\index{s@$\sigmak$} is defined by \(\sigmak(B)\defg\int_{\ra}t\sigmah(\dif t)\) and belongs to \(\Mggqaag{\ra}{\seq{\sau{j}}{j}{0}{\kappa-1}}\) (see \eqref{s+} and \eqref{so}).
\erema

\bremnl{see~\zitaa{MR3324594}{\clemp{4.3}{11}}}{R1017}
 Let \(\kappa\in\NOinf \) and let \(\seqska \in\Kgqkappa\), then
 \[
  \ophua{n}{0}
  =(-1)^n\rk*{\rprod_{k=0}^{n-1}\DSpm{k}\DSpl{k}}^\inv
 \]
 for all \(n\in\N\) with \(2n-1\leq\kappa\), and
 \[
  \spkua{n}{0}
  =(-1)^n\ek*{\rk*{\rprod_{k=0}^{n-1}\DSpm{k}\DSpl{k}}\DSpm{n}}^\inv
 \]
 for all \(n\in\N\) with \(2n\leq\kappa\).
\erema

We have the following representation of \(\Uu{m}\) in terms of the
above introduced matrix polynomials:
\bpropnl{see~\zitaa{MR3324594}{\cthmp{4.7}{68}}}{P1422}
 Let \(\kappa\in\NOinf \) and let \(\seqska \in\Kgqkappa\). For all \(z\in\C\), then
 \[
  \Uua{2n}{z}
  =
  \bMat
   \spkua{n}{z}&-\sphua{n}{z}\\
   -z\opkua{n}{z}&\ophua{n}{z}
  \eMat
  \bMat
   [\spkua{n}{0}]^\inv&\Oqq\\
   \Oqq&[\ophua{n}{0}]^\inv
  \eMat
 \]
 for all \(n\in\NO\) with \(2n\leq\kappa\), and
 \[
  \Uua{2n+1}{z}
  =
  \bMat
   \spkua{n}{z}&-\sphua{n+1}{z}\\
   -z\opkua{n}{z}&\ophua{n+1}{z}
  \eMat
  \bMat
   [\spkua{n}{0}]^\inv&\Oqq\\
   \Oqq&[\ophua{n+1}{0}]^\inv
  \eMat
 \]
 for all \(n\in\NO\) with \(2n+1\leq\kappa\).
\eprop

\section{The \hKt{}}\label{S1649}
 In~\zitaa{114arxiv}{\S7} a transformation for sequences of complex \tpqa{matrices} was considered using the following concept of \tsraute{s} presented in~\zita{MR3014197}. The paper~\zita{MR3014197} deals with the question of invertibility as it applies to matrix sequences.

\begin{defn}\label{D1430}
 Let \(\kappa \in \NOinf \) and let \(\seqska \) be a sequence of complex \tpqa{matrices}. The sequence \(\seq{\su{j}^\rez}{j}{0}{\kappa}\)\index{$\seq{\su{j}^\rez}{j}{0}{\kappa}$} given recursively by
 \begin{align*}
  s_0^\rez&\defg\su{0}^\MP&
  &\text{and}&
  s_j^\rez&\defg-\su{0}^\MP \sum_{\ell=0}^{j-1} s_{j-\ell} \su{\ell}^\rez
 \end{align*}
 for all \(j\in\mn{1}{\kappa}\) is called the \noti{\tsrautea{\seqska }}{r@\tsraute}.
\end{defn}

\begin{defn}\label{D1059}
 Let \(\kappa\in\Ninf \) and let \(\seqska \) be a sequence of complex \tpqa{matrices}. Then, the sequence \(\seq{\su{j}^\Sta{1}}{j}{0}{\kappa-1}\)\index{$\seq{\su{j}^\Sta{1}}{j}{0}{\kappa-1}$} given by
 \[
  \su{j}^\Sta{1}
  \defg-\su{0}\su{j+1}^\rez\su{0}
 \]
 is called the \noti{first \tKta{\(\seqska \)}}{s@\hKt{}!first}.%
\end{defn}
Observe that this transformation coincides with the first \tlasnt{\alpha} from~\zitaa{114arxiv}{\cdefnp{7.1}{33}} for \(\alpha=0\) and served  together with its counterpart for matrix-valued functions in~\zitas{114arxiv,141arxiv} as elementary step of a Schur type algorithm to solve the truncated matricial moment problem on the semi-infinite interval \([\alpha,\infp)\).
\bdefil{D1632}
 Let \(\kappa\in\NOinf \) and let \(\seqska \) be a sequence of complex \tpqa{matrices}. The sequence \(\seq{\su{j}^\Sta{0}}{j}{0}{\kappa}\) given by \(\su{j}^\Sta{0}\defg\su{j}\) for all \(j\in\mn{0}{\kappa}\) is called the \noti{\taKta{0}{\(\seqska \)}}{s@\hKt{}!$0$th}. In the case \(\kappa\geq1\), for all \(k\in\mn{1}{\kappa}\), the \noti{\taKt{k} \(\seq{\su{j}^\Sta{k}}{j}{0}{\kappa-k}\) of \(\seqska \)}{s@\hKt{}!$k$th}\index{$\seq{\su{j}^\Sta{k}}{j}{0}{\kappa-k}$} is recursively defined by
 \[
  \su{j}^\Sta{k}
  \defg t_j^\Sta{1}
 \]
 for all \(j\in\mn{0}{\kappa-k}\), where \(\seq{t_j}{j}{0}{\kappa-(k-1)}\) denotes the \taKta{(k-1)}{\(\seqska \)}.
\edefi

\bpropnl{see~\zitaa{114arxiv}{\cthmp{8.10(c)}{42}}}{P1542c}
 Let \(\kappa\in\NOinf \) and \(\seqska \in\Kgqu{\kappa}\). Then \(\seq{\su{j}^\Sta{k}}{j}{0}{\kappa-k}\in\Kgqu{\kappa-k}\) for all \(k\in\mn{0}{\kappa}\).
\eprop

\bpropnl{see~\zitaa{114arxiv}{\cthmp{8.10(d)}{42}}}{P1542d}
 Let \(m\in\NO\) and \(\seqsinf\in\Kggdoq{m}\). Then \(\seq{\su{j}^\Sta{k}}{j}{0}{\kappa-k}\in\Kggdoq{\max\set{0,m-k}}\) for all \(k\in\mn{0}{\kappa}\).
\eprop

\bpropnl{see~\zitaa{114arxiv}{\cthmp{8.10(e)}{42}}}{P1542e}
 Let \(\kappa\in\NOinf \) and \(\seqska \in\Kggdqu{\kappa}\). Then \(\seq{\su{j}^\Sta{k}}{j}{0}{\kappa-k}\in\Kggdqu{\kappa-k}\) for all \(k\in\mn{0}{\kappa}\).
\eprop

\bpropl{T1658}
 Let \(\kappa\in\NOinf \), let \(\seqska \in\Kgqu{\kappa}\) with \tSp{} \(\seq{\Spu{j}}{j}{0}{\kappa}\), and let \(k\in\mn{0}{\kappa}\). Then \(\seq{\Spu{k+j}}{j}{0}{\kappa-k}\) is the \tSpa{\(\seq{\su{j}^\Sta{k}}{j}{0}{\kappa-k}\)}.
\eprop
\bproof
 According to~\zitaa{MR3014201}{\cpropp{2.20}{221}} we have \(\seqska \in\Kggeqkappa\). Hence, the application of~\zitaa{114arxiv}{\cthmp{9.26}{57}} yields the assertion.
\eproof

\section{The \hdsp{} after \hKt{ation}}\label{S1011}
 From now on we consider only infinite sequences \(\seqsinf\in\Kgqinf\). Let \(\ell\in\NO\). According to \rprop{P1542c}, then the \taKt{\ell} \(\seq{\su{j}^\Sta{\ell}}{j}{0}{\infi}\) of \(\seqsinf\) belongs to \(\Kgqinf\). If \(X\) is an object build from the sequence \(\seqsinf\), then we will use the notation \(X^{(\ell)}\)\index{\(X^{(\ell)}\)} for this object build from the sequence \(\seq{\su{j}^\Sta{\ell}}{j}{0}{\infi}\), e.\,g., in view of \eqref{HK}, we have\index{h@$\Hu{n}^\Sta{\ell}\defg\matauuo{\su{j+k}^\Sta{\ell}}{j,k}{0}{n}$}\index{k@$\Ku{n}^\Sta{\ell}\defg\matauuo{\su{j+k+1}^\Sta{\ell}}{j,k}{0}{n}$}
 \begin{align*}
  \Hu{n}^{(\ell)}&\defg\matauuo{\su{j+k}^\Sta{\ell}}{j,k}{0}{n}&
  &\text{and}&
  \Ku{n}^{(\ell)}&\defg\matauuo{\su{j+k+1}^\Sta{\ell}}{j,k}{0}{n}
 \end{align*}
 for all \(n\in\NO\). From \rdefi{D1021} and \rprop{T1658} we see:
\bremal{R2352}
 Let \(\seqsinf\in\Kgqinf\). Then \(\Lu{k}^{(1)}=\Lau{k}\) and \(\Lau{k}^{(1)}=\Lu{k+1}\) for all \(k\in\NO\).
\erema

 In view of \rprop{P1542c} and \rdefi{D1455}, we are particularly able to introduce the following notation:
\bnotal{N1045}
 Let \(\seqsinf \in\Kgqinf\) and let \(\ell\in\NO\). Then, we write \symb{[\seq{\DSpol{k}{\ell}}{k}{0}{\infi},\seq{\DSpom{k}{\ell}}{k}{0}{\infi}]} for the \tdspa{the \taKta{\ell}{\(\seqsinf\)}}.
\enota

 From \rtheo{T1513} and \rprop{T1658} we obtain then:
\blemml{L1656}
 Let \(\seqsinf \in\Kgqinf\) and let \(\ell\in\NO\), then
 \begin{align*}
  \DSpol{k}{\ell}
  &=\rk*{\rprod_{j=0}^k\Spu{2j+\ell}\Spu{2j+\ell+1}^\inv}\Spu{2k+\ell+1}\rk*{\rprod_{j=0}^k\Spu{2j+\ell}\Spu{2j+\ell+1}^\inv}^\ad
 \shortintertext{and}
  \DSpom{k}{\ell}
  &=
  \begin{cases}
    \Spu{\ell}^\inv\incase{k=0}\\
    \rk*{\rprod_{j=0}^{k-1}\Spu{2j+\ell}^\inv\Spu{2j+\ell+1}}\Spu{2k+\ell}^\inv\rk*{\rprod_{j=0}^{k-1}\Spu{2j+\ell}^\inv\Spu{2j+\ell+1}}^\ad\incase{k\geq1}
  \end{cases}
 \end{align*}
 for all \(k\in\NO\).
\elemm

Using \rlemm{L1656} we conclude:
\blemml{L0759}
 Let \(\seqsinf \in\Kgqinf\) and let \(\ell\in\NO\), then
 \begin{align*}
  \DSpol{k}{\ell}
  &=\rk*{\rprod_{j=0}^{m-1}\Spu{2j+\ell}\Spu{2j+\ell+1}^\inv}\DSpol{k-m}{\ell+2m}\rk*{\rprod_{j=0}^{m-1}\Spu{2j+\ell}\Spu{2j+\ell+1}^\inv}^\ad
  \shortintertext{and}
  \DSpom{k}{\ell}
  &=\rk*{\rprod_{j=0}^{m-1}\Spu{2j+\ell}^\inv\Spu{2j+\ell+1}}\DSpom{k-m}{\ell+2m}\rk*{\rprod_{j=0}^{m-1}\Spu{2j+\ell}\Spu{2j+\ell+1}^\inv}^\ad
 \end{align*}
 for all \(k\in\N\) and \(m\in\mn{1}{k}\).
\elemm

\blemml{L0855}
 Let \(\seqsinf \in\Kgqinf\) and let \(\ell\in\NO\), then
 \begin{align*}
  \DSpol{k}{\ell+1}&=\Spu{\ell}\DSpom{k+1}{\ell}\Spu{\ell}&
  &\text{and}&
  \DSpom{k}{\ell+1}&=\Spu{\ell}^\inv\DSpol{k}{\ell}\Spu{\ell}^\inv
 \end{align*}
 for all \(k\in\NO\).
\elemm
\bproof
 From \rprop{T1337} we know that the matrices \(\Spu{j}\) are \tH{} and invertible for all \(j\in\NO\). Using \rlemm{L1656}, we obtain thus
 \[
  \begin{split}
   \DSpol{k}{\ell+1}
   &=\Spu{\ell}\Spu{\ell}^\inv\DSpol{k}{\ell+1}\Spu{\ell}^\invad\Spu{\ell}^\ad\\
   &=\Spu{\ell}\Spu{\ell}^\inv\rk*{\rprod_{j=0}^k\Spu{2j+\ell+1}\Spu{2j+\ell+2}^\inv}\Spu{2k+\ell+2}\rk*{\rprod_{j=0}^k\Spu{2j+\ell+1}\Spu{2j+\ell+2}^\inv}^\ad\Spu{\ell}^\invad\Spu{\ell}^\ad\\
   &=\Spu{\ell}\rk*{\rprod_{j=0}^k\Spu{2j+\ell}^\inv\Spu{2j+\ell+1}}\Spu{2k+\ell+2}^\inv\Spu{2k+\ell+2}\Spu{2k+\ell+2}^\invad\rk*{\rprod_{j=0}^k\Spu{2j+\ell}^\inv\Spu{2j+\ell+1}}^\ad\Spu{\ell}^\ad\\
   &=\Spu{\ell}\rk*{\rprod_{j=0}^k\Spu{2j+\ell}^\inv\Spu{2j+\ell+1}}\Spu{2(k+1)+\ell}^\inv\rk*{\rprod_{j=0}^k\Spu{2j+\ell}^\inv\Spu{2j+\ell+1}}^\ad\Spu{\ell}
   =\Spu{\ell}\DSpom{k+1}{\ell}\Spu{\ell}
  \end{split}
 \]
 for all \(k\in\NO\) and
 \[
  \begin{split}
   \DSpom{0}{\ell+1}
   =\Spu{\ell}^\inv\Spu{\ell}\DSpom{0}{\ell+1}\Spu{\ell}^\ad\Spu{\ell}^\invad
   &=\Spu{\ell}^\inv\Spu{\ell}\Spu{\ell+1}^\inv\Spu{\ell}^\ad\Spu{\ell}^\invad
   =\Spu{\ell}^\inv\Spu{\ell}\Spu{\ell+1}^\inv\Spu{\ell+1}\Spu{\ell+1}^\invad\Spu{\ell}^\ad\Spu{\ell}^\invad\\
   &=\Spu{\ell}^\inv(\Spu{\ell}\Spu{\ell+1}^\inv)\Spu{\ell+1}(\Spu{\ell}\Spu{\ell+1}^\inv)^\ad\Spu{\ell}^\inv
   =\Spu{\ell}^\inv\DSpol{0}{\ell}\Spu{\ell}^\inv
  \end{split}
 \]
 and furthermore
 \[
  \begin{split}
   &\DSpom{k}{\ell+1}
   =\Spu{\ell}^\inv\Spu{\ell}\DSpom{k}{\ell+1}\Spu{\ell}^\ad\Spu{\ell}^\invad\\
   &=\Spu{\ell}^\inv\Spu{\ell}\rk*{\rprod_{j=0}^{k-1}\Spu{2j+\ell+1}^\inv\Spu{2j+\ell+2}}\Spu{2k+\ell+1}^\inv\rk*{\rprod_{j=0}^{k-1}\Spu{2j+\ell+1}^\inv\Spu{2j+\ell+2}}^\ad\Spu{\ell}^\ad\Spu{\ell}^\invad\\
   &=\Spu{\ell}^\inv\Spu{\ell}\rk*{\rprod_{j=0}^{k-1}\Spu{2j+\ell+1}^\inv\Spu{2j+\ell+2}}\Spu{2k+\ell+1}^\inv\Spu{2k+\ell+1}\Spu{2k+\ell+1}^\invad\rk*{\rprod_{j=0}^{k-1}\Spu{2j+\ell+1}^\inv\Spu{2j+\ell+2}}^\ad\Spu{\ell}^\ad\Spu{\ell}^\inv\\
   &=\Spu{\ell}^\inv\rk*{\rprod_{j=0}^k\Spu{2j+\ell}\Spu{2j+\ell+1}^\inv}\Spu{2k+\ell+1}\rk*{\rprod_{j=0}^k\Spu{2j+\ell}\Spu{2j+\ell+1}^\inv}^\ad\Spu{\ell}^\inv
   =\Spu{\ell}^\inv\DSpol{k}{\ell}\Spu{\ell}^\inv
  \end{split}
 \]
 for all \(k\in\N\).
\eproof

In view of \rdefi{D1021} and \eqref{L}, we obtain from \rlemm{L0855} the following relation between the \tdsp{s} of a \tSpd{} sequence and its first \tKt{}:
\btheol{R0921}
 Let \(\seqsinf\in\Kgqinf\) with \tdsp{} \([\seq{\DSpl{k}}{k}{0}{\infty},\seq{\DSpm{k}}{k}{0}{\infty}]\). Then, the \tdsp{} \([\seq{\DSpol{k}{1}}{k}{0}{\infty},\seq{\DSpom{k}{1}}{k}{0}{\infty}]\) of \(\seq{\su{j}^\Sta{1}}{j}{0}{\infty}\) is given by
 \begin{align*}
  \DSpol{k}{1}&=\su{0}\DSpm{k+1}\su{0}&
  &\text{and}&
  \DSpom{k}{1}&=\su{0}^\inv\DSpl{k}\su{0}^\inv
 \end{align*}
 for all \(k\in\NO\).
\etheo

\section{The resolvent matrix corresponding to a \hKt{ed} moment sequence}\label{S0834}
 In view of \rtheo{T1523} and~\zitaa{MR3133464}{\cremp{8.4}{24}} we have:
\bremal{R1054}
 Let \(\seqsinf \in\Kgqinf\) with \tdsp{} \([\seq{\DSpl{k}}{k}{0}{\infty},\seq{\DSpm{k}}{k}{0}{\infty}]\) and \tSp{} \(\seq{\Spu{j}}{j}{0}{\infi}\), then
 \begin{align*}
  \rprod_{j=0}^m\Spu{2j}^\inv\Spu{2j+1}&=\rk*{\rprod_{j=0}^m\DSpm{j}\DSpl{j}}^\inv&
  &\text{and}&
  \rprod_{j=0}^m\Spu{2j}\Spu{2j+1}^\inv&=\rk*{\rprod_{j=0}^m\DSpm{j}\DSpl{j}}^\ad
 \end{align*}
 for all \(m\in\NO\).
\erema

\blemml{L0926}
 Let \(\seqsinf \in\Kgqinf\) with \tdsp{} \([\seq{\DSpl{k}}{k}{0}{\infty},\seq{\DSpm{k}}{k}{0}{\infty}]\) and \tStiqosoomp{} \(\PQPQ\). Then,
 \[
  \DSpl{k}
  =\ek*{\ophua{m}{0}}^\invad\DSpol{k-m}{2m}\ek*{\ophua{m}{0}}^\inv
  =\spkua{m}{0}\DSpom{k-m}{2m+1}\ek*{\spkua{m}{0}}^\ad
 \]
 and
 \begin{align*}
  \DSpm{k}&=\ophua{m}{0}\DSpom{k-m}{2m}\ek*{\ophua{m}{0}}^\ad,&
  \DSpm{k+1}=\ek*{\spkua{m}{0}}^\invad\DSpol{k-m}{2m+1}\ek*{\spkua{m}{0}}^\inv
 \end{align*}
 for all \(k\in\NO\) and \(m\in\mn{0}{k}\).
\elemm
\bproof
 Let \(k\in\NO\) and \(m\in\mn{0}{k}\). In the case \(m=0\) all assertions follow from \rnota{N1104} and \rrema{R0921}. Now suppose that \(k\in\N\) and \(m\in\mn{1}{k}\). From \rremass{R1054}{R1017} we obtain
 \begin{align*}
  \rprod_{j=0}^{m-1}\Spu{2j}^\inv\Spu{2j+1}
  &=(-1)^m\ophua{m}{0}&
  &\text{and}&
  \rprod_{j=0}^{m-1}\Spu{2j}\Spu{2j+1}^\inv&=(-1)^m\ek*{\ophua{m}{0}}^\invad.
 \end{align*}
 Hence, using \rlemm{L0759} with \(\ell=0\), we get
 \[\begin{split}
  \DSpl{k}
  &=\rk*{\rprod_{j=0}^{m-1}\Spu{2j}\Spu{2j+1}^\inv}\DSpol{k-m}{2m}\rk*{\rprod_{j=0}^{m-1}\Spu{2j}\Spu{2j+1}^\inv}^\ad
  =\ek*{\ophua{m}{0}}^\invad\DSpol{k-m}{2m}\ek*{\ophua{m}{0}}^\inv
 \end{split}\]
 and
 \begin{equation}\label{L0926.1}\begin{split}
  \DSpm{k}
  &=\rk*{\rprod_{j=0}^{m-1}\Spu{2j}^\inv\Spu{2j+1}}\DSpom{k-m}{2m}\rk*{\rprod_{j=0}^{m-1}\Spu{2j}\Spu{2j+1}^\inv}^\ad
  =\ophua{m}{0}\DSpom{k-m}{2m}\ek*{\ophua{m}{0}}^\ad.
 \end{split}\end{equation}
 According to \rlemm{L0855} we have
 \begin{align*}
  \DSpol{k-m}{2m+1}&=\Spu{2m}\DSpom{k-m+1}{2m}\Spu{2m}&
  &\text{and}&
  \DSpom{k-m}{2m+1}&=\Spu{2m}^\inv\DSpol{k-m}{2m}\Spu{2m}^\inv.
 \end{align*}
 \rtheo{T1523} and \rrema{R1017} yield
 \[
  \Spu{2m}
  =\rk*{\rprod_{j=0}^{m-1}\DSpm{j}\DSpl{j}}^\invad\DSpm{m}^\inv\rk*{\rprod_{j=0}^{m-1}\DSpm{j}\DSpl{j}}^\inv
  =\ek*{\ophua{m}{0}}^\ad\spkua{m}{0}.
 \]
 Thus, we get
 \[
  \begin{split}
   \ek*{\ophua{m}{0}}^\invad\DSpol{k-m}{2m}\ek*{\ophua{m}{0}}^\inv
   &=\ek*{\ophua{m}{0}}^\invad\Spu{2m}\DSpom{k-m}{2m+1}\Spu{2m}^\ad\ek*{\ophua{m}{0}}^\inv\\
   &=\spkua{m}{0}\DSpom{k-m}{2m+1}\ek*{\spkua{m}{0}}^\ad
  \end{split}
 \]
 and, using the already proved identity \eqref{L0926.1} with \(k+1\) instead of \(k\), furthermore
 \[
   \begin{split}
    \DSpm{k+1}
    =\ophua{m}{0}\DSpom{k-m+1}{2m}\ek*{\ophua{m}{0}}^\ad
    &=\ophua{m}{0}\Spu{2m}^\invad\DSpol{k-m}{2m+1}\Spu{2m}^\inv\ek*{\ophua{m}{0}}^\ad\\
    &=\ek*{\spkua{m}{0}}^\invad\DSpol{k-m}{2m+1}\ek*{\spkua{m}{0}}^\inv.\qedhere
   \end{split}
 \]
\eproof

Let \(\seqsinf \in\Kgqinf\). Then, for all \(n\in\NO\), the matrices \(\ophua{n}{0}\) and \(\spkua{n}{0}\)
are invertible according to \rrema{R1017} and \rnota{N1104}. Thus, for all \(n\in\NO\), let
\beql{PP}
 \PPu{n}
 \defg
 \bMat
  [\ophua{n}{0}]^\invad&\Oqq\\
  \Oqq&\ophua{n}{0}
 \eMat
\eeq
\index{p@$\PPu{n}$}and let \(\QQu{n}\colon\C\to\Coo{2q}{2q}\)\index{q@$\QQu{n}$} be defined by
\beql{QQ}
 \QQua{n}{z}
 \defg
 \bMat
  \Oqq&\spkua{n}{0}\\
  -z[\spkua{n}{0}]^\invad&\Oqq
 \eMat.
\eeq
Obviously, the matrix \(\PPu{n}\) is invertible for all
\(n\in\NO\) and \(\QQua{n}{z}\) is invertible for all \(n\in\NO\)
and all \(z\in\C\setminus\set{0}\). Furthermore, let
\begin{align}\label{LLMMl}
 \Lluo{k}{\ell}
 &\defg
 \bMat
    \Iq&\DSpol{k}{\ell}\\
    \Oqq&\Iq
   \eMat&
&\text{and}&
 \Mmuo{k}{\ell}
 &\defg
 \bMat
   \Iq&\Oqq\\
   -z\DSpom{k}{\ell}&\Iq
  \eMat
\end{align}
for all \(k,\ell\in\NO\) and all \(z\in\C\), where \([\seq{\DSpol{k}{\ell}}{k}{0}{\infi},\seq{\DSpom{k}{\ell}}{k}{0}{\infi}]\) denotes the \tdspa{the \taKta{\ell}{\(\seqsinf\)}}.

\blemml{L1402}
 Let \(\seqsinf \in\Kgqinf\), then
 \[
  \Llu{k}
  =\PPu{m}\Lluo{k-m}{2m}\PPu{m}^\inv
  =\QQu{m}\Mmuo{k-m}{2m+1}\QQu{m}^\inv
 \]
 and
 \begin{align*}
  \Mmu{k}&=\PPu{m}\Mmuo{k-m}{2m}\PPu{m}^\inv,&
  \Mmu{k+1}=\QQu{m}\Lluo{k-m}{2m+1}\QQu{m}^\inv
 \end{align*}
 for all \(k\in\NO\) and \(m\in\mn{0}{k}\).
\elemm
\bproof
 Let \(k\in\NO\), let \(m\in\mn{0}{k}\), and let \(z\in\C\setminus\set{0}\). In view of \rlemm{L0926} we have then
 \[
  \begin{split}
   &\PPu{m}\Lluo{k-m}{2m}\PPu{m}^\inv\\
   &=
   \bMat
    [\ophua{m}{0}]^\invad&\Oqq\\
    \Oqq&\ophua{m}{0}
   \eMat
   \bMat
    \Iq&\DSpol{k-m}{2m}\\
    \Oqq&\Iq
   \eMat
   \bMat
    [\ophua{m}{0}]^\invad&\Oqq\\
    \Oqq&\ophua{m}{0}
   \eMat^\inv\\
   &=
   \bMat
    [\ophua{m}{0}]^\invad&[\ophua{m}{0}]^\invad\DSpol{k-m}{2m}\\
    \Oqq&\ophua{m}{0}
   \eMat
   \bMat
    [\ophua{m}{0}]^\ad&\Oqq\\
    \Oqq&[\ophua{m}{0}]^\inv
   \eMat\\
   &=
   \bMat
    \Iq&[\ophua{m}{0}]^\invad\DSpol{k-m}{2m}[\ophua{m}{0}]^\inv\\
    \Oqq&\Iq
   \eMat
   =
   \bMat
    \Iq&\DSpl{k}\\
    \Oqq&\Iq
   \eMat
   =\Llu{k}
  \end{split}
 \]
 and
 \[
  \begin{split}
   &\QQua{m}{z}\Mmuoa{k-m}{2m+1}{z}\ek*{\QQua{m}{z}}^\inv\\
   &=
   \bMat
    \Oqq&\spkua{m}{0}\\
    -z[\spkua{m}{0}]^\invad&\Oqq
   \eMat
   \bMat
    \Iq&\Oqq\\
    -z\DSpom{k-m}{2m+1}&\Iq
   \eMat
   \bMat
    \Oqq&\spkua{m}{0}\\
    -z[\spkua{m}{0}]^\invad&\Oqq
   \eMat^\inv\\
   &=
   \bMat
    -z\spkua{m}{0}\DSpom{k-m}{2m+1}&\spkua{m}{0}\\
    -z[\spkua{m}{0}]^\invad&\Oqq
   \eMat
   \bMat
    \Oqq&-z^\inv[\spkua{m}{0}]^\ad\\
    [\spkua{m}{0}]^\inv&\Oqq
   \eMat\\
   &=
   \bMat
    \Iq&\spkua{m}{0}\DSpom{k-m}{2m+1}[\spkua{m}{0}]^\ad\\
    \Oqq&\Iq
   \eMat
   =
   \bMat
    \Iq&\DSpl{k}\\
    \Oqq&\Iq
   \eMat
   =\Llu{k}
  \end{split}
 \]
 and furthermore
 \[
  \begin{split}
   &\PPu{m}\Mmuoa{k-m}{2m}{z}\PPu{m}^\inv\\
   &=
   \bMat
    [\ophua{m}{0}]^\invad&\Oqq\\
    \Oqq&\ophua{m}{0}
   \eMat
   \bMat
    \Iq&\Oqq\\
    -z\DSpom{k-m}{2m}&\Iq
   \eMat
   \bMat
    [\ophua{m}{0}]^\invad&\Oqq\\
    \Oqq&\ophua{m}{0}
   \eMat^\inv\\
   &=
   \bMat
    [\ophua{m}{0}]^\invad&\Oqq\\
    -z\ophua{m}{0}\DSpom{k-m}{2m}&\ophua{m}{0}
   \eMat
   \bMat
    [\ophua{m}{0}]^\ad&\Oqq\\
    \Oqq&[\ophua{m}{0}]^\inv
   \eMat\\
   &=
   \bMat
    \Iq&\Oqq\\
    -z\ophua{m}{0}\DSpom{k-m}{2m}[\ophua{m}{0}]^\ad&\Iq
   \eMat
   =
   \bMat
    \Iq&\Oqq\\
    -z\DSpm{k}&\Iq
   \eMat
   =\Mmua{k}{z}
  \end{split}
 \]
 and
 \begin{equation*}
  \begin{split}
   &\QQua{m}{z}\Lluo{k-m}{2m+1}\ek*{\QQua{m}{z}}^\inv\\
   &=
   \bMat
    \Oqq&\spkua{m}{0}\\
    -z[\spkua{m}{0}]^\invad&\Oqq
   \eMat
   \bMat
    \Iq&\DSpol{k-m}{2m+1}\\
    \Oqq&\Iq
   \eMat
   \bMat
    \Oqq&\spkua{m}{0}\\
    -z[\spkua{m}{0}]^\invad&\Oqq
   \eMat^\inv\\
   &=
   \bMat
    \Oqq&\spkua{m}{0}\\
    -z[\spkua{m}{0}]^\invad&-z[\spkua{m}{0}]^\invad\DSpol{k-m}{2m+1}
   \eMat
   \bMat
    \Oqq&-z^\inv[\spkua{m}{0}]^\ad\\
    [\spkua{m}{0}]^\inv&\Oqq
   \eMat\\
   &=
   \bMat
    \Iq&\Oqq\\
    -z[\spkua{m}{0}]^\invad\DSpol{k-m}{2m+1}[\spkua{m}{0}]^\inv&\Iq
   \eMat
   =
   \bMat
    \Iq&\Oqq\\
    -z\DSpm{k+1}&\Iq
   \eMat=\Mmua{k+1}{z}.\qedhere
  \end{split}
 \end{equation*}
\eproof

 In view of \rprop{P1542c} and \rdefi{D1137}, we are able to introduce the following notations:
\bnotal{N1053}
 Let \(\seqsinf \in\Kgqinf\) and let \(\ell\in\NO\). Then, we write \symb{[\seq{\ophuo{k}{\ell}}{k}{0}{\infi},\seq{\sphuo{k}{\ell}}{k}{0}{\infi},\seq{\opkuo{k}{\ell}}{k}{0}{\infi},\seq{\spkuo{k}{\ell}}{k}{0}{\infi}]} for the \tStiqosoompo{the \taKta{\ell}{\(\seqsinf\)}}.
\enota

 Let \(\seqsinf \in\Kgqinf\) and \(\ell\in\NO\). Then, the matrices \(\ophuoa{n}{\ell}{0}\) and \(\spkuoa{n}{\ell}{0}\) are invertible. In accordance with \eqref{PP} and \eqref{QQ}, let
\beql{PPell}
 \PPuo{n}{\ell}
 \defg
 \bMat
  [\ophuoa{n}{\ell}{0}]^\invad&\Oqq\\
  \Oqq&\ophuoa{n}{\ell}{0}
 \eMat
\eeq
\index{p@$\PPuo{n}{\ell}$}for all \(n\in\NO\), and let \(\QQuo{n}{\ell}\colon\C\to\Coo{2q}{2q}\)\index{q@$\QQuo{n}{\ell}$} be defined by
\beql{QQell}
 \QQuoa{n}{\ell}{z}
 \defg
 \bMat
  \Oqq&\spkuoa{n}{\ell}{0}\\
  -z[\spkuoa{n}{\ell}{0}]^\invad&\Oqq
 \eMat.
\eeq
 Obviously, the matrix \(\PPuo{n}{\ell}\) is invertible for all \(n\in\NO\) and \(\QQuoa{n}{\ell}{z}\) is invertible for all \(n\in\NO\) and all \(z\in\C\setminus\set{0}\). For all \(m\in\NO\), let\index{q@$\Qq{m}$}
\beql{Qq}
 \Qq{m}
 \defg\rprod_{\ell=0}^m\QQuo{0}{\ell}.
\eeq
 In particular, we have \(\Qq{0}=\QQuo{0}{0}=\QQu{0}\).

 In view of \rprop{P1542c}, we can apply \rlemm{L1402} with \(m=0\) to the \taKta{\ell}{a sequence \(\seqsinf\in\Kgqinf\)} and obtain:
\bremal{R1422}
 Let \(\seqsinf\in\Kgqinf\), then
 \begin{align*}
  \Lluo{k}{\ell}\QQuo{0}{\ell}&=\QQuo{0}{\ell}\Mmuo{k}{\ell+1}&
 &\text{and}&
  \Mmuo{k+1}{\ell}\QQuo{0}{\ell}&=\QQuo{0}{\ell}\Lluo{k}{\ell+1}
 \end{align*}
 for all \(k,\ell\in\NO\).
\erema

 By repeated application of \rrema{R1422} we get:
\blemml{L1446}
 Let \(\seqsinf\in\Kgqinf\), then
 \begin{align*}
  \Qq{2n}\Lluo{k}{2n+1}&=\Mmu{k+n+1}\Qq{2n},&
  \Qq{2n}\Mmuo{k}{2n+1}&=\Llu{k+n}\Qq{2n}
 \shortintertext{and}
  \Qq{2n+1}\Lluo{k}{2n+2}&=\Llu{k+n+1}\Qq{2n+1},&
  \Qq{2n+1}\Mmuo{k}{2n+2}&=\Mmu{k+n+1}\Qq{2n+1}
 \end{align*}
 for all \(k,n\in\NO\).
\elemm
\bproof
 Using \rrema{R1422} and \eqref{Qq}, we obtain
 \begin{align*}
  \Qq{0}\Lluo{k}{1}
  &=\QQuo{0}{0}\Lluo{k}{0+1}
  =\Mmuo{k+1}{0}\QQuo{0}{0}
  =\Mmu{k+1}\Qq{0},\\
  \Qq{0}\Mmuo{k}{1}
  &=\QQuo{0}{0}\Mmuo{k}{0+1}
  =\Lluo{k}{0}\QQuo{0}{0}
  =\Llu{k}\Qq{0},\\
  \Qq{1}\Lluo{k}{2}
  &=\QQuo{0}{0}\QQuo{0}{1}\Lluo{k}{1+1}\\
  &=\QQuo{0}{0}\Mmuo{k+1}{1}\QQuo{0}{1}
  =\QQuo{0}{0}\Mmuo{k+1}{0+1}\QQuo{0}{1}
  =\Lluo{k+1}{0}\QQuo{0}{0}\QQuo{0}{1}
  =\Llu{k+1}\Qq{1},
 \shortintertext{and}
  \Qq{1}\Mmuo{k}{2}
  &=\QQuo{0}{0}\QQuo{0}{1}\Mmuo{k}{1+1}\\
  &=\QQuo{0}{0}\Lluo{k}{1}\QQuo{0}{1}
  =\QQuo{0}{0}\Lluo{k}{0+1}\QQuo{0}{1}
  =\Mmuo{k+1}{0}\QQuo{0}{0}\QQuo{0}{1}
  =\Mmu{k+1}\Qq{1}.
 \end{align*}
 for all \(k\in\NO\). Now let \(n\in\N\) and suppose that
 \begin{align*}
  \Qq{2n-1}\Lluo{j}{2n}&=\Llu{j+n}\Qq{2n-1},&
  \Qq{2n-1}\Mmuo{j}{2n}&=\Mmu{j+n}\Qq{2n-1}
 \end{align*}
 hold true for all \(j\in\NO\). Taking additionally into account \eqref{Qq} and \rrema{R1422}, we get for all \(k\in\NO\) then
 \[\begin{split}
  \Qq{2n}\Lluo{k}{2n+1}
  =\Qq{2n-1}\QQuo{0}{2n}\Lluo{k}{2n+1}
  &=\Qq{2n-1}\Mmuo{k+1}{2n}\QQuo{0}{2n}\\
  &=\Mmu{k+1+n}\Qq{2n-1}\QQuo{0}{2n}
  =\Mmu{k+n+1}\Qq{2n}
 \end{split}\]
 and
 \[\begin{split}
  \Qq{2n}\Mmuo{k}{2n+1}
  =\Qq{2n-1}\QQuo{0}{2n}\Mmuo{k}{2n+1}
  &=\Qq{2n-1}\Lluo{k}{2n}\QQuo{0}{2n}\\
  &=\Llu{k+n}\Qq{2n-1}\QQuo{0}{2n}
  =\Llu{k+n}\Qq{2n}.
 \end{split}\]
 Using this and again \eqref{Qq} and \rrema{R1422}, we obtain for all \(\ell\in\NO\) furthermore
 \[\begin{split}
  \Qq{2n+1}\Lluo{\ell}{2n+2}
  &=\Qq{2n}\QQuo{0}{2n+1}\Lluo{\ell}{(2n+1)+1}\\
  &=\Qq{2n}\Mmuo{\ell+1}{2n+1}\QQuo{0}{2n+1}
  =\Llu{\ell+1+n}\Qq{2n}\QQuo{0}{2n+1}
  =\Llu{\ell+n+1}\Qq{2n+1}
 \end{split}\]
 and
 \[\begin{split}
  \Qq{2n+1}\Mmuo{\ell}{2n+2}
  &=\Qq{2n}\QQuo{0}{2n+1}\Mmuo{\ell}{(2n+1)+1}\\
  &=\Qq{2n}\Lluo{\ell}{2n+1}\QQuo{0}{2n+1}
  =\Mmu{\ell+n+1}\Qq{2n}\QQuo{0}{2n+1}
  =\Mmu{\ell+n+1}\Qq{2n+1}.
 \end{split}\]
 Thus, the assertion is proved by mathematical induction.
\eproof

\bnotal{N1053AA}
 Let \(\seqsinf \in\Kgqinf\). For all \(\ell\in\NO\) denote by \(\seq{\Uuo{m}{\ell}}{m}{0}{\infi}\) the \tsqDmpo{the \taKta{\ell}{\(\seqsinf\)}}. Then, for all \(k,\ell\in\NO\), let \symba{\UCR{k}{\ell}\defg\Uuo{k-\ell}{\ell}}{u}.
\enota

 In particular, this means \(\UCR{k}{0}\defg\Uu{k}\) and, in view of \eqref{Uu2n} and \eqref{Uu2n1}, furthermore
 \begin{align}
  \UCR{2n}{2k}
  &\defg
  \bMat
   \pauo{n-k}{2k}&\pbuo{n-k}{2k}\\
   \pcuo{n-k}{2k}&\pduo{n-k}{2k}\\
  \eMat,&
  \UCR{2n}{2k+1}
  &\defg
  \bMat
   \pauo{n-k-1}{2k+1}&\pbuo{n-k}{2k+1}\\
   \pcuo{n-k-1}{2k+1}&\pduo{n-k}{2k+1}\\
  \eMat,\label{aaa66a}
 \shortintertext{and}
  \UCR{2n+1}{2k}
  &\defg
  \bMat
   \pauo{n-k}{2k}&\pbuo{n-k+1}{2k}\\
   \pcuo{n-k}{2k}&\pduo{n-k+1}{2k}\\
  \eMat,&
  \UCR{2n+1}{2k+1}
  &\defg
  \bMat
   \pauo{n-k}{2k+1}&\pbuo{n-k}{2k+1}\\
   \pcuo{n-k}{2k+1}&\pduo{n-k}{2k+1}\\
  \eMat\label{aaa77a}
 \end{align}
 for all \(n,k\in\NO\).
%
%

 Now we are able to write down the connection between the resolvent matrices \(\UCR{m}{0}=\Uu{m}\) and \(\UCR{m}{1}=\Uuo{m-1}{1}\) given via  \eqref{aaa66a} and \eqref{aaa77a}.

\bpropl{P1443}
 Let \(\seqsinf\in\Kgqinf\). Then
\(
 \UCR{m}{1}
 =\QQu{0}^\inv\Mmu{0}^\inv\UCR{m}{0}\QQu{0}
\)
 for all \(m\in\N\).
\eprop
\bproof
 According to \rprop{P1542c}, we have \(\seq{\su{j}^\Sta{1}}{j}{0}{\infi}\in\Kgqinf\). \rprop{P1033} yields
 \begin{align*}
  \UCR{1}{0}&=\Mmu{0}\Llu{0}&
  &\text{and}&
  \UCR{1}{1}&=\Mmuo{0}{1}.
 \end{align*}
 Hence, using \rlemm{L1402}, we obtain
 \[\begin{split}
  \QQu{0}^\inv\Mmu{0}^\inv\UCR{1}{0}\QQu{0}
  &=\QQu{0}^\inv\Mmu{0}^\inv\Mmu{0}\Llu{0}\QQu{0}\\
  &=\QQu{0}^\inv\Mmu{0}^\inv\Mmu{0}\QQu{0}\Mmuo{0}{1}\QQu{0}^\inv\QQu{0}
  =\Mmuo{0}{1}
  =\UCR{1}{1}.
 \end{split}\]
 Let \(n\in\N\). \rprop{P1033} yields then
 \begin{align*}
  \UCR{2n}{0}&=\rk*{\rprod_{k=0}^{n-1}\Mmu{k}\Llu{k}}\Mmu{n},&
  \UCR{2n+1}{0}&=\rprod_{k=0}^{n}\rk{\Mmu{k}\Llu{k}}
 \shortintertext{and}
  \UCR{2n}{1}&=\rprod_{k=0}^{n-1}(\Mmuo{k}{1}\Lluo{k}{1}),&
  \UCR{2n+1}{1}&=\rk*{\rprod_{k=0}^{n-1}\Mmuo{k}{1}\Lluo{k}{1}}\Mmuo{n}{1}.
 \end{align*}
 Hence, using \rlemm{L1402}, we obtain
 \[
  \begin{split}
  \QQu{0}^\inv\Mmu{0}^\inv\UCR{2n}{0}\QQu{0}
  &=\QQu{0}^\inv\Mmu{0}^\inv\rk*{\rprod_{k=0}^{n-1}\Mmu{k}\Llu{k}}\Mmu{n}\QQu{0}
  =\QQu{0}^\inv\Mmu{0}^\inv\Mmu{0}\rk*{\rprod_{k=0}^{n-1}\Llu{k}\Mmu{k+1}}\QQu{0}\\
  &=\QQu{0}^\inv\rk*{\rprod_{k=0}^{n-1}\QQu{0}\Mmuo{k}{1}\QQu{0}^\inv\QQu{0}\Lluo{k}{1}\QQu{0}^\inv}\QQu{0}
  =\rprod_{k=0}^{n-1}(\Mmuo{k}{1}\Lluo{k}{1})
  =\UCR{2n}{1}
  \end{split}
 \]
 and
 \[
  \begin{split}
  \QQu{0}^\inv\Mmu{0}^\inv\UCR{2n+1}{0}\QQu{0}
  &=\QQu{0}^\inv\Mmu{0}^\inv\rk*{\rprod_{k=0}^{n}\Mmu{k}\Llu{k}}\QQu{0}
  =\QQu{0}^\inv\Mmu{0}^\inv\Mmu{0}\rk*{\rprod_{k=0}^{n-1}\Llu{k}\Mmu{k+1}}\Llu{n}\QQu{0}\\
  &=\QQu{0}^\inv\rk*{\rprod_{k=0}^{n-1}\QQu{0}\Mmuo{k}{1}\QQu{0}^\inv\QQu{0}\Lluo{k}{1}\QQu{0}^\inv}\QQu{0}\Mmuo{n}{1}\QQu{0}^\inv\QQu{0}\\
  &=\rk*{\rprod_{k=0}^{n-1}\Mmuo{k}{1}\Lluo{k}{1}}\Mmuo{n}{1}
  =\UCR{2n+1}{1}.\qedhere
  \end{split}
 \]
\eproof

 In view of \rprop{P1542c}, we can apply \rprop{P1443} to the \taKta{\ell}{a sequence \(\seqsinf\in\Kgqinf\)} and obtain:
\begin{rem}\label{aaaa400}
 Let \(\seqsinf\in\Kgqinf\). Then
\(
 \UCR{m}{\ell}
 =\Mmuo{0}{\ell}\QQuo{0}{\ell}\UCR{m}{\ell+1}\rk{\QQuo{0}{\ell}}^\inv
\)
 for all \(m\in\N\) and all \(\ell\in\mn{0}{m-1}\).
\end{rem}


 By applying \rrema{aaaa400} twice, we obtain:
\blemml{cor00aa00}
 Let \(\seqsinf\in\Kgqinf\). Then
\[
 \UCR{m}{0}
 =\Mmu{0}\Llu{0}\rk{\QQu{0}\QQuo{0}{1}}\UCR{m}{2}\rk{\QQu{0}\QQuo{0}{1}}^\inv
\]
 for all \(m\in\minf{2}\).
\elemm
\bproof
 According to \rprop{P1542c} the sequences \(\seq{\su{j}^\Sta{1}}{j}{0}{\infi}\) and \(\seq{\su{j}^\Sta{2}}{j}{0}{\infi}\) both belong to \(\Kgqinf\). The application of \rrema{aaaa400} to \(\seqsinf\) and \(\seq{\su{j}^\Sta{1}}{j}{0}{\infi}\) yields
 \begin{align*}
  \UCR{m}{0}&=\Mmu{0}\QQu{0}\UCR{m}{1}\QQu{0}^\inv&
  &\text{and}&
  \UCR{m}{1}&=\Mmuo{0}{1}\QQuo{0}{1}\UCR{m}{2}\rk{\QQuo{0}{1}}^\inv.
 \end{align*}
 Hence, \(\UCR{m}{0}=\Mmu{0}\QQu{0}\Mmuo{0}{1}\QQuo{0}{1}\UCR{m}{2}\rk{\QQuo{0}{1}}^\inv\QQu{0}^\inv\). From \rrema{R1422} with \(k=\ell=0\), we obtain furthermore \(\QQu{0}\Mmuo{0}{1}=\Llu{0}\QQu{0}\) which completes the proof.
\eproof

%

\blemml{L1105}
 Let \(\seqsinf \in\Kgqinf\), then
 \[
  \UCR{m}{0}\Qq{m}
  =\rprod_{\ell=0}^m\rk{\Mmuo{0}{\ell}\QQuo{0}{\ell}}
 \]
 for all \(m\in\NO\).
\elemm
\bproof
 According to \rprop{P1542c}, we have \(\seq{\su{j}^\Sta{\ell}}{j}{0}{\infi}\in\Kgqinf\) for all \(\ell\in\NO\). \rprop{P1033} yields
 \begin{align*}
  \Uu{0}&=\Mmu{0}&
  &\text{and}&
  \Uu{1}&=\Mmu{0}\Llu{0}.
 \end{align*}
 Hence, we obtain
 \[
  \Uu{0}\Qq{0}
  =\Mmu{0}\QQu{0}
  =\Mmuo{0}{0}\QQuo{0}{0}
 \]
 and, using \rrema{R1422}, furthermore
 \[
  \Uu{1}\Qq{1}
  =\Mmu{0}\Llu{0}\QQuo{0}{0}\QQuo{0}{1}
  =\Mmu{0}\QQu{0}\Mmuo{0}{1}\QQuo{0}{1}
  =\rk{\Mmuo{0}{0}\QQuo{0}{0}}\rk{\Mmuo{0}{1}\QQuo{0}{1}}.
 \]
 Now let \(n\in\N\) and suppose that
 \[
  \Uu{2n-1}\Qq{2n-1}
  =\rprod_{\ell=0}^{2n-1}\rk{\Mmuo{0}{\ell}\QQuo{0}{\ell}}
 \]
 holds true. \rprop{P1033} yields
 \begin{align*}
  \Uu{2n}&=\Uu{2n-1}\Mmu{n}&
 &\text{and}&
  \Uu{2n+1}&=\Uu{2n}\Llu{n}.
 \end{align*}
 Hence, using \rlemm{L1446} with \(k=0\), we obtain
 \[\begin{split}
  \Uu{2n}\Qq{2n}
  =\Uu{2n-1}\Mmu{n}\Qq{2n-1}\QQuo{0}{2n}
  =\Uu{2n-1}\Qq{2n-1}\Mmuo{0}{2n}\QQuo{0}{2n}
  =\rprod_{\ell=0}^{2n}\rk{\Mmuo{0}{\ell}\QQuo{0}{\ell}}
  \end{split}
 \]
 and with that furthermore
 \[
  \Uu{2n+1}\Qq{2n+1}
  =\Uu{2n}\Llu{n}\Qq{2n}\QQuo{0}{2n+1}
  =\Uu{2n}\Qq{2n}\Mmuo{0}{2n+1}\QQuo{0}{2n+1}
  =\rprod_{\ell=0}^{2n+1}\rk{\Mmuo{0}{\ell}\QQuo{0}{\ell}}.
 \]
 In view of \(\Uu{m}=\UCR{m}{0}\) for all \(m\in\NO\), the assertion is thus proved by mathematical induction.
\eproof

 The following result expresses the resolvent matrix $\UCR{m}{0}=\Uu{m}$ in terms of the resolvent matrices \(\UCR{\ell}{0}=\Uu{\ell}\) and \(\UCR{m}{\ell+1}=\Uuo{m-\ell-1}{\ell+1}\), which can be considered as a splitting of the original moment problem:
\btheol{Taaa007}
  Let \(\seqsinf\in\Kgqinf\). Then
\(
 \UCR{m}{0}
 =\UCR{\ell}{0}\Qq{\ell}\UCR{m}{\ell+1}\Qq{\ell}^\inv
\)
 for all \(m\in\N\) and \(\ell\in\mn{0}{m-1}\).
\etheo
\bproof
 Let \(m\in\N\). According to \rprop{P1033} and \rrema{aaaa400} we have
 \[
  \UCR{0}{0}\Qq{0}\UCR{m}{1}\Qq{0}^\inv
  =\Uuo{0}{0}\QQuo{0}{0}\UCR{m}{1}\rk{\QQuo{0}{0}}^\inv
  =\Uu{0}\QQu{0}\UCR{m}{1}\QQu{0}^\inv
  =\Mmu{0}\QQu{0}\UCR{m}{1}\QQu{0}^\inv
  =\UCR{m}{0}.
 \]
 Now suppose \(m\geq2\) and that
 \[
  \UCR{m}{0}
  =\UCR{\ell}{0}\Qq{\ell}\UCR{m}{\ell+1}\Qq{\ell}^\inv
 \]
 holds true for some \(\ell\in\mn{0}{m-2}\). Using \rlemm{L1105} and \rrema{aaaa400}, we can conclude then
 \[\begin{split}
  \UCR{\ell+1}{0}\Qq{\ell+1}\UCR{m}{\ell+2}\Qq{\ell+1}^\inv
  &=\ek*{\rprod_{k=0}^{\ell+1}\rk{\Mmuo{0}{k}\QQuo{0}{k}}}\UCR{m}{\ell+2}\rk{\Qq{\ell}\QQuo{0}{\ell+1}}^\inv\\
  &=\ek*{\rprod_{k=0}^{\ell}\rk{\Mmuo{0}{k}\QQuo{0}{k}}}\rk{\Mmuo{0}{\ell+1}\QQuo{0}{\ell+1}}\UCR{m}{(\ell+1)+1}\rk{\QQuo{0}{\ell+1}}^\inv\Qq{\ell}^\inv\\
  &=\UCR{\ell}{0}\Qq{\ell}\UCR{m}{\ell+1}\Qq{\ell}^\inv
  =\UCR{m}{0}.
 \end{split}\]
 Thus, the assertion is proved by mathematical induction.
\eproof

\section[Orthogonal matrix polynomials corresponding to a transformed sequence]{Orthogonal matrix polynomials corresponding to a \hKt{ed} moment sequence}\label{S1945}
%

 From \rpropss{P1443}{P1422} we obtain the following connection between the matrix polynomials \(\ophu{n}\), \(\sphu{n}\), \(\opku{n}\), and \(\spku{n}\), and the matrix polynomials \(\ophuo{n}{1}\), \(\sphuo{n}{1}\), \(\opkuo{n}{1}\), and \(\spkuo{n}{1}\) given in \rnotass{N1104}{N1053} and \rdefi{D1137} corresponding to a \tSpd{} sequence and its first \tKt, respectively:

\bpropl{C1536}
 Let \(\seqsinf\in\Kgqinf\) with first \tKt{} $(s_j^\Sta{1})_{j=0}^\infi$. Then the \tStiqosoomp{}  \(\PQPQ\) of \(\seqsinf\) and \([\seq{\ophuo{k}{1}}{k}{0}{\infi},\seq{\sphuo{k}{1}}{k}{0}{\infi},\seq{\opkuo{k}{1}}{k}{0}{\infi},\seq{\spkuo{k}{1}}{k}{0}{\infi}]\) of $(s_j^\Sta{1})_{j=0}^\infi$, resp., fulfill
 \begin{align*}
  \ek*{\ophuoa{n}{1}{z}}\ek*{\ophuoa{n}{1}{0}}^\inv&=\ek*{\su{0}^\inv\spkua{n}{z}}\ek*{\su{0}^\inv\spkua{n}{0}}^\inv,\\
  \ek*{\sphuoa{n}{1}{z}}\ek*{\ophuoa{n}{1}{0}}^\inv&=\ek*{\spkua{n}{z}-\su{0}\opkua{n}{z}}\ek*{\su{0}^\inv\spkua{n}{0}}^\inv,\\
  \ek*{\opkuoa{n-1}{1}{z}}\ek*{\spkuoa{n-1}{1}{0}}^\inv&=\ek*{\su{0}^\inv\sphua{n}{z}}\ek*{-\su{0}\ophua{n}{0}}^\inv,
 \shortintertext{and}
  \ek*{\spkuoa{n-1}{1}{z}}\ek*{\spkuoa{n-1}{1}{0}}^\inv&
  =\ek*{z\sphua{n}{z}-\su{0}\ophua{n}{z}}\ek*{-\su{0}\ophua{n}{0}}^\inv.
 \end{align*}
 for all \(n\in\N\) and all \(z\in\C\).
\eprop
\bproof
 According to \rprop{P1542c}, we have \(\seq{\su{j}^\Sta{1}}{j}{0}{\infty}\in\Kgqinf\).
  In view of \eqref{ma11}, \eqref{QQ}, and \(\su{0}^\ad=\su{0}\), we get
 \begin{align*}
  \Mmua{0}{z}
  &=
  \bMat
   \Iq&\Oqq\\
   -z\su{0}^\inv&\Iq
  \eMat&
 &\text{and}&
  \QQua{0}{z}
  &=
  \bMat
   \Oqq&\su{0}\\
   -z\su{0}^\inv&\Oqq
  \eMat.
 \end{align*}
 Hence,
 \[
  \Mmua{0}{z}\QQua{0}{z}
  =
  \bMat
   \Oqq&\su{0}\\
   -z\su{0}^\inv&-z\Iq
  \eMat
 \]
 and thus
 \[
  \ek*{\QQua{0}{z}}^\inv\ek*{\Mmua{0}{z}}^\inv
  =\ek*{\Mmua{0}{z}\QQua{0}{z}}^\inv
  =
  \bMat
   -\Iq&-z^\inv\su{0}\\
   \su{0}^\inv&\Oqq
  \eMat.
 \]
 Using \eqref{aaa66a} and \rprop{P1443} we obtain then
 \[\begin{split}
  \bMat
   \pauoa{n-1}{1}{z}&\pbuoa{n}{1}{z}\\
   \pcuoa{n-1}{1}{z}&\pduoa{n}{1}{z}\\
  \eMat
  &=\UCRa{2n}{1}{z}
   =\ek*{\QQua{0}{z}}^\inv\ek*{\Mmua{0}{z}}^\inv\UCRa{2n}{0}{z}\QQua{0}{z}\\
   &=
   \bMat
    -\Iq&-z^\inv\su{0}\\
    \su{0}^\inv&\Oqq
   \eMat
   \bMat
    \pauoa{n}{0}{z}&\pbuoa{n}{0}{z}\\
    \pcuoa{n}{0}{z}&\pduoa{n}{0}{z}\\
   \eMat
   \bMat
    \Oqq&\su{0}\\
    -z\su{0}^\inv&\Oqq
   \eMat\\
   &=
   \bMat
    \ek{z\pbuoa{n}{0}{z}+\su{0}\pduoa{n}{0}{z}}\su{0}^\inv&-\ek{\pauoa{n}{0}{z}+z^\inv\su{0}\pcuoa{n}{0}{z}}\su{0}\\
    -z\su{0}^\inv\pbuoa{n}{0}{z}\su{0}^\inv&\su{0}^\inv\pauoa{n}{0}{z}\su{0}
   \eMat.
  \end{split}
 \]
 Since \eqref{aaa66a} and \rprop{P1422} yield
 \[
  \bMat
   \pauoa{n}{0}{z}&\pbuoa{n}{0}{z}\\
   \pcuoa{n}{0}{z}&\pduoa{n}{0}{z}\\
  \eMat
  =
  \bMat
   \ek{\spkua{n}{z}}[\spkua{n}{0}]^\inv&-\ek{\sphua{n}{z}}[\ophua{n}{0}]^\inv\\
   -z\ek{\opkua{n}{z}}[\spkua{n}{0}]^\inv&\ek{\ophua{n}{z}}[\ophua{n}{0}]^\inv
  \eMat
 \]
 and
 \[
  \bMat
   \pauoa{n-1}{1}{z}&\pbuoa{n}{1}{z}\\
   \pcuoa{n-1}{1}{z}&\pduoa{n}{1}{z}\\
  \eMat
  =
  \bMat
   \ek{\spkuoa{n-1}{1}{z}}[\spkuoa{n-1}{1}{0}]^\inv&-\ek{\sphuoa{n}{1}{z}}[\ophuoa{n}{1}{0}]^\inv\\
   -z\ek{\opkuoa{n-1}{1}{z}}[\spkuoa{n-1}{1}{0}]^\inv&\ek{\ophuoa{n}{1}{z}}[\ophuoa{n}{1}{0}]^\inv
  \eMat
 \]
 we can conclude
 \begin{align*}
  &\ek*{\spkuoa{n-1}{1}{z}}\ek*{\spkuoa{n-1}{1}{0}}^\inv
  =\ek*{z\rk*{-\ek*{\sphua{n}{z}}\ek*{\ophua{n}{0}}^\inv}+\su{0}\ek*{\ophua{n}{z}}\ek*{\ophua{n}{0}}^\inv}\su{0}^\inv,\\
  -&\ek*{\sphuoa{n}{1}{z}}\ek*{\ophuoa{n}{1}{0}}^\inv
  =-\ek*{\ek*{\spkua{n}{z}}\ek*{\spkua{n}{0}}^\inv+z^\inv\su{0}\rk*{-z\ek*{\opkua{n}{z}}\ek*{\spkua{n}{0}}^\inv}}\su{0},\\
  -z&\ek*{\opkuoa{n-1}{1}{z}}\ek*{\spkuoa{n-1}{1}{0}}^\inv
  =-z\su{0}^\inv\rk*{-\ek*{\sphua{n}{z}}\ek*{\ophua{n}{0}}^\inv}\su{0}^\inv,
\shortintertext{and}
  &\ek*{\ophuoa{n}{1}{z}}\ek*{\ophuoa{n}{1}{0}}^\inv
  =\su{0}^\inv\rk*{\ek*{\spkua{n}{z}}\ek*{\spkua{n}{0}}^\inv}\su{0}.
 \end{align*}
 Hence, the asserted identities follow.
\eproof

\blemml{L0822}
 Let \(\seqsinf\in\Kgqinf\). Then
 \begin{align*}
  \ophuoa{n}{1}{0}&=\su{0}^\inv\spkua{n}{0}&
 &\text{and}&
  \spkuoa{n}{1}{0}&=-\su{0}\ophua{n+1}{0}
 \end{align*}
 for all \(n\in\N\).
\elemm
\bproof
 According to \rprop{P1542c}, we have \(\seq{\su{j}^\Sta{1}}{j}{0}{\infty}\in\Kgqinf\). Denote by \([\seq{\DSpl{k}}{k}{0}{\infty},\seq{\DSpm{k}}{k}{0}{\infty}]\) and \([\seq{\DSpol{k}{1}}{k}{0}{\infty},\seq{\DSpom{k}{1}}{k}{0}{\infty}]\) the \tdspa{\(\seqsinf\) and \(\seq{\su{j}^\Sta{1}}{j}{0}{\infty}\)}, respectively. Let \(n\in\N\). Using \rremass{R1017}{R0921} we have, in view of \eqref{M0}, then
 \[
  \begin{split}
   \ophuoa{n}{1}{0}
   &=(-1)^n\rk*{\rprod_{k=0}^{n-1}\DSpom{k}{1}\DSpol{k}{1}}^\inv
   =(-1)^n\rk*{\rprod_{k=0}^{n-1}\su{0}^\inv\DSpl{k}\su{0}^\inv\su{0}\DSpm{k+1}\su{0}}^\inv\\
   &=(-1)^n\ek*{\su{0}^\inv\DSpl{0}\rk*{\rprod_{\ell=1}^{n-1}\DSpm{\ell}\DSpl{\ell}}\DSpm{n}\su{0}}^\inv
   =(-1)^n\ek*{\DSpm{0}\DSpl{0}\rk*{\rprod_{\ell=1}^{n-1}\DSpm{\ell}\DSpl{\ell}}\DSpm{n}\su{0}}^\inv\\
   &=\su{0}^\inv\rk*{(-1)^n\ek*{\rk*{\rprod_{\ell=0}^{n-1}\DSpm{\ell}\DSpl{\ell}}\DSpm{n}}^\inv}
   =\su{0}^\inv\spkua{n}{0}
  \end{split}
 \]
 and
 \[
  \begin{split}
   \spkuoa{n}{1}{0}
   &=(-1)^n\ek*{\rk*{\rprod_{k=0}^{n-1}\DSpom{k}{1}\DSpol{k}{1}}\DSpom{n}{1}}^\inv\\
   &=(-1)^n\ek*{\rk*{\rprod_{k=0}^{n-1}\su{0}^\inv\DSpl{k}\su{0}^\inv\su{0}\DSpm{k+1}\su{0}}\su{0}^\inv\DSpl{n}\su{0}^\inv}^\inv\\
   &=(-1)^n\ek*{\su{0}^\inv\DSpl{0}\rk*{\rprod_{\ell=1}^{n}\DSpm{\ell}\DSpl{\ell}}\su{0}^\inv}^\inv
   =(-1)^n\ek*{\DSpm{0}\DSpl{0}\rk*{\rprod_{\ell=1}^{n}\DSpm{\ell}\DSpl{\ell}}\su{0}^\inv}^\inv\\
   &=-\su{0}\ek*{(-1)^{n+1}\rk*{\rprod_{\ell=0}^{n}\DSpm{\ell}\DSpl{\ell}}^\inv}
   =-\su{0}\ophua{n+1}{0}.\qedhere
  \end{split}
 \]
\eproof

\btheol{T0820}
 Let \(\seqsinf\in\Kgqinf\) with first \tKt{} $(s_j^\Sta{1})_{j=0}^\infi$. Then the \tStiqosoomp{} \(\PQPQ\) of \(\seqsinf\) and \([\seq{\ophuo{k}{1}}{k}{0}{\infi},\seq{\sphuo{k}{1}}{k}{0}{\infi},\seq{\opkuo{k}{1}}{k}{0}{\infi},\seq{\spkuo{k}{1}}{k}{0}{\infi}]\) of $(s_j^\Sta{1})_{j=0}^\infi$, resp., fulfill
 \begin{align*}
  \ophuoa{n}{1}{z}&=\su{0}^\inv\spkua{n}{z},&\sphuoa{n}{1}{z}&=\spkua{n}{z}-\su{0}\opkua{n}{z},\\
  \opkuoa{n-1}{1}{z}&=\su{0}^\inv\sphua{n}{z},&\spkuoa{n-1}{1}{z}&=z\sphua{n}{z}-\su{0}\ophua{n}{z}
 \end{align*}
 for all \(n\in\N\) and all \(z\in\C\).
\etheo
\bproof
 Combine \rprop{C1536} with \rlemm{L0822}.
\eproof

\bpropl{orthg-001}
 Let \(\seqsinf\in\Kgqinf\) with \tStiqosoomp{}  \(\PQPQ\) and first \tKt{} $(s_j^\Sta{1})_{j=0}^\infi$. Then the sequences $(s_j^\Sta{1})_{j=0}^\infi$ and $(s_{j+1}^\Sta{1})_{j=0}^\infi$ both belong to \(\Kgqinf\). Furthermore \(\seq{\su{0}^\inv\spku{k}}{k}{0}{\infi}\) is a \tmrosmpa{$(s_j^\Sta{1})_{j=0}^\infi$} and \(\seq{\su{0}^\inv\sphu{k+1}}{k}{0}{\infi}\) is a \tmrosmpa{$(s_{j+1}^\Sta{1})_{j=0}^\infi$}. In particular, if \(\muh\in\Mggqaag{\ra}{(s_j^\Sta{1})_{j=0}^\infi}\), we have the orthogonality relations
 \begin{align*}
  \int_{\ra}\ek*{\su{0}^\inv\spkua{m}{t}}^\ad\muh(\dif t)\ek*{\su{0}^\inv\spkua{n}{t}}
  &=
  \begin{cases}
   \Oqq\incase{m\neq n}\\
   \Lau{n}\incase{m=n}
  \end{cases}
 \shortintertext{and}
  \int_{\ra}\ek*{\su{0}^\inv\sphua{m}{t}}^\ad\muk(\dif t)\ek*{\su{0}^\inv\sphua{n}{t}}
  &=
  \begin{cases}
   \Oqq\incase{m\neq n}\\
   \Lu{n}\incase{m=n}
  \end{cases}
 \end{align*}
 for all \(m,n\in\NO\), where \(\muk\colon\Bra\to\Cggq\) is defined by \symba{\muk(B)\defg\int_{\ra}t\muh(\dif t)}{m} and belongs to \(\Mggqaag{\ra}{(s_{j+1}^\Sta{1})_{j=0}^\infi}\).
\eprop
\bproof
 In view of \rprop{P1542c} and \rrema{R2347}, we see that $(s_j^\Sta{1})_{j=0}^\infi$ and $(s_{j+1}^\Sta{1})_{j=0}^\infi$ both belong to \(\Kgqinf\). Since \(\spkua{0}{z}=\su{0}\) for all \(z\in\C\), the combination of \rtheo{T0820} with \rremass{R0001}{R2352} completes the proof.
\eproof

 Let \(m\in\NO\) and let \(\seqs{m}\in\Kgqu{m}\). Then we want to draw the attention to two distinguished elements of the solution set \(\Mggqaakg{\ra}{\seq{\su{j}}{j}{0}{m}}\). This concerns those measures \(\smin{m}\)\index{s@\(\smin{m}\)} and \(\smax{m}\)\index{s@\(\smax{m}\)}, respectively, the \tSt{s} of which are generated via \rthmp{T1327}{T1327.a} by the two constant pairs \((\iota,\theta),(\theta,\iota)\in\SPq\), where \(\iota\) and \(\theta\) are the constant functions in \(\Cs\) with values \(\Iq\) and \(\Oqq\), respectively. The measures \(\smin{m}\) and \(\smax{m}\) are called the \noti{lower}{extremal element of $\Mggqaakg{\ra}{\seq{\su{j}}{j}{0}{m}}$!lower} and \noti{upper extremal elements of \(\Mggqaakg{\ra}{\seq{\su{j}}{j}{0}{m}}\)}{extremal element of $\Mggqaakg{\ra}{\seq{\su{j}}{j}{0}{m}}$!upper}, respectively. In view of \rtheo{T1327}, \eqref{Uu2n}, \eqref{Uu2n1}, and \rprop{P1422}, for \(n\in\NO\) and \(z\in\Cs\) we infer for the corresponding \tSt{s}
\begin{align}
 \sttrhla{\smin{2n}}{z}&=\ek*{\paua{n}{z}}\ek*{\pcua{n}{z}}^\inv=-\ek*{\spkua{n}{z}}\ek*{z\opkua{n}{z}}^\inv,\label{min0}\\
 \sttrhla{\smax{2n}}{z}&=\ek*{\pbua{n}{z}}\ek*{\pdua{n}{z}}^\inv=-\ek*{\sphua{n}{z}}\ek*{\ophua{n}{z}}^\inv\label{max0}
\shortintertext{and}
 \sttrhla{\smin{2n+1}}{z}&=\ek*{\paua{n}{z}}\ek*{\pcua{n}{z}}^\inv=-\ek*{\spkua{n}{z}}\ek*{z\opkua{n}{z}}^\inv,\label{min1}\\
 \sttrhla{\smax{2n+1}}{z}&=\ek*{\pbua{n+1}{z}}\ek*{\pdua{n+1}{z}}^\inv=-\ek*{\sphua{n+1}{z}}\ek*{\ophua{n+1}{z}}^\inv.\label{max1}
\end{align}
 The functions introduced in \eqref{min0}--\eqref{max1} play an important role in the considerations of Yu.~M.~Dyukarev~\zita{MR2053150}. We refer the reader to~\zitaa{MR2053150}{\cSect{3}} for a detailed discussion of these functions and their extremality properties.

 Taking into account \rnota{N1104}, we see that \(\sttrhl{\smin{2n+1}}\) coincides with \(\sttrhl{\smin{2n}}\) and does not depend on \(\su{2n+1}\) and that \(\sttrhl{\smax{2n}}\) coincides with \(\sttrhl{\smax{2n-1}}\) and does not depend on \(\su{2n}\). In particular
 \begin{align}\label{o=e}
  \smin{2n+1}&=\smin{2n}&
 &\text{and}&
  \smax{2n}&=\smax{2n-1}.
 \end{align}
 Since \(\smin{m}\) and \(\smax{m}\) both belong to $\Mggqaakg{\ra}{\seq{\su{j}}{j}{0}{m}}$, we can hence conclude
 \begin{align}\label{sm=}
  \smin{2n}&\in\MggqAAg{\ra}{\seq{\su{j}}{j}{0}{2n}},&
  \smax{2n-1}&\in\MggqAAg{\ra}{\seq{\su{j}}{j}{0}{2n-1}}.
 \end{align}
 The  lower and upper extremal elements \(\smin{m}\) and \(\smax{m}\) of \(\Mggqaakg{\ra}{\seq{\su{j}}{j}{0}{m}}\) are concentrated on a finite number of points in \(\ra\). In particular, they possess power moments up to any order, which coincide with \tzext{s} introduced in \rdefi{D0909}:

\blemml{L1014}
 Let \(n\in\N\) and let \(\seqs{2n-1}\in\Kgqu{2n-1}\). Then \(\int_\ra x^j\smax{2n-1}(\dif x)=\zext{s}{j}\) for all \(j\in\NO\), where \(\seqzinf\) is the \tzexto{\(\seqs{2n-1}\)}.
\elemm
\bproof
 For all \(j\in\NO\) let \(t_j\defg\int_{\ra}x^j\smax{2n-1}(\dif x)\). According to \rtheo{T1121}, we have then \(\seq{t_j}{j}{0}{\infi}\in\Kggqinf\). Because of \eqref{sm=}, we get furthermore \(\su{j}=t_j\) for all \(j\in\mn{0}{2n-1}\). From \rprop{P0926} and \eqref{L} we can consequently conclude that the matrix \(t_{2n}-\Theta_{2n}\) is \tnnH{}, where \(\Theta_{2n}\defg\zuu{n}{2n-1}\Hu{n-1}^\inv\yuu{n}{2n-1}\). Now, we consider an arbitrary \(\epsilon>0\). Let \(\su{2n}\defg\Theta_{2n}+\epsilon\Iq\) and denote by \((\Spu{j})_{j=0}^{2n}\) the \tSpa{\(\seqs{2n}\)}. In view of \eqref{L}, the matrix \(\Spu{2n}=\Lu{n}=\epsilon\Iq\) is \tpH{}. From \rprop{T1337} we then can easily conclude \(\seqs{2n}\in\Kgqu{2n}\). As an element of \(\Mggqaakg{\ra}{\seqs{2n}}\), the measure \(\smax{2n}\) fulfills \(\su{2n}-\int_\ra x^{2n}\smax{2n}(\dif x)\in\Cggq\). Using \eqref{o=e}, we obtain
 \[
  \Theta_{2n}
  \leq t_{2n}
  =\int_\ra x^{2n}\smax{2n}(\dif x)
  \leq\su{2n}
  =\Theta_{2n}+\epsilon\Iq.
 \]
 Since this holds true for all \(\epsilon>0\), we get \(t_{2n}=\Theta_{2n}\). Thus, the sequence \(\seq{t_j}{j}{0}{\infi}\) belongs to \(\Kggdoq{2n}\). Hence, we can apply \rlemm{L0912} to see that \(\seq{t_j}{j}{0}{\infi}\) is the \tzexto{\(\seq{t_j}{j}{0}{2n-1}\)}. Because of \(\su{j}=t_j\) for all \(j\in\mn{0}{2n-1}\), then \(t_j=\zext{s}{j}\) for all \(j\in\NO\), which completes the proof.
\eproof

\blemml{L1037}
 Let \(n\in\NO\) and let \(\seqs{2n}\in\Kgqu{2n}\). Then \(\int_\ra x^j\smin{2n}(\dif x)=\zext{s}{j}\) for all \(j\in\NO\), where \(\seqzinf\) is the \tzexto{\(\seqs{2n}\)}.
\elemm
\bproof
 For all \(j\in\NO\) let \(t_j\defg\int_{\ra}x^j\smin{2n}(\dif x)\). According to \rtheo{T1121}, we have then \(\seq{t_j}{j}{0}{\infi}\in\Kggqinf\). Because of \eqref{sm=}, we get furthermore \(\su{j}=t_j\) for all \(j\in\mn{0}{2n}\). From \rprop{P0926} and \eqref{La} we can consequently conclude that the matrix \(t_{2n+1}-\Theta_{2n+1}\) is \tnnH{}, where \(\Theta_1\defg\Oqq\) and \(\Theta_{2n+1}\defg\zuu{n+1}{2n}\Ku{n-1}^\inv\yuu{n+1}{2n}\) for \(n\in\N\). Now, we consider an arbitrary \(\epsilon>0\). Let \(\su{2n+1}\defg\Theta_{2n+1}+\epsilon\Iq\) and denote by \((\Spu{j})_{j=0}^{2n+1}\) the \tSpa{\(\seqs{2n+1}\)}. In view of \eqref{La}, the matrix \(\Spu{2n+1}=\Lau{n}=\epsilon\Iq\) is \tpH{}. From \rprop{T1337} we then can easily conclude \(\seqs{2n+1}\in\Kgqu{2n+1}\). As an element of \(\Mggqaakg{\ra}{\seqs{2n+1}}\), the measure \(\smin{2n+1}\) fulfills \(\su{2n+1}-\int_\ra x^{2n+1}\smin{2n+1}(\dif x)\in\Cggq\). Using \eqref{o=e}, we obtain
 \[
  \Theta_{2n+1}
  \leq t_{2n+1}
  =\int_\ra x^{2n+1}\smin{2n+1}(\dif x)
  \leq\su{2n+1}
  =\Theta_{2n+1}+\epsilon\Iq.
 \]
 Since this holds true for all \(\epsilon>0\), we get \(t_{2n+1}=\Theta_{2n+1}\). Thus, the sequence \(\seq{t_j}{j}{0}{\infi}\) belongs to \(\Kggdoq{2n+1}\). Hence, we can apply \rlemm{L0912} to see that \(\seq{t_j}{j}{0}{\infi}\) is the \tzexto{\(\seq{t_j}{j}{0}{2n}\)}. Because of \(\su{j}=t_j\) for all \(j\in\mn{0}{2n}\), then \(t_j=\zext{s}{j}\) for all \(j\in\NO\), which completes the proof.
\eproof

 In combination with \eqref{min0}--\eqref{max1}, \rtheo{T0820} yields a relation between the lower and upper extremal elements associated with a \tSpd{} sequence and its first \tKt{}:

\bpropl{P1058}
 Let \(\seqsinf\in\Kgqinf\) with first \tKt{} $(s_j^\Sta{1})_{j=0}^\infi$. Then \(\seqs{m}\) and $(s_j^\Sta{1})_{j=0}^m$ both belong to \(\Kgqu{m}\) for all \(m\in\NO\). For all \(m\in\NO\) denote by \(\smin{m}\) and \(\smax{m}\) the lower and upper extremal element of \(\Mggqaakg{\ra}{\seq{\su{j}}{j}{0}{m}}\). Furthermore, for all \(m\in\N\), let \(\smino{m}{1}\) and \(\smaxo{m}{1}\) be the lower and upper extremal element of \(\Mggqaakg{\ra}{\seq{\su{j}^\Sta{1}}{j}{0}{m}}\). Then
\begin{align}
 \sttrhla{\smino{2n-2}{1}}{z}&=\sttrhla{\smino{2n-1}{1}}{z}=-\su{0}-\su{0}\ek*{z\sttrhla{\smax{2n}}{z}}^\inv\su{0}\label{P1058.B1}
 \shortintertext{and}
 \sttrhla{\smaxo{2n-1}{1}}{z}&=\sttrhla{\smaxo{2n}{1}}{z}=-\su{0}-\su{0}\ek*{z\sttrhla{\smin{2n}}{z}}^\inv\su{0}\label{P1058.B2}
\end{align}
 for all \(n\in\N\) and all \(z\in\Cs\).
\eprop
\bproof
 According to \rprop{P1542c}, we have \(\seq{\su{j}^\Sta{1}}{j}{0}{\infty}\in\Kgqinf\). Denote by \(\PQPQ\) the \tStiqosoomp{} of \(\seqsinf\) and by \([\seq{\ophuo{k}{1}}{k}{0}{\infi},\seq{\sphuo{k}{1}}{k}{0}{\infi},\seq{\opkuo{k}{1}}{k}{0}{\infi},\seq{\spkuo{k}{1}}{k}{0}{\infi}]\) the \tStiqosoomp{} of $(s_j^\Sta{1})_{j=0}^\infi$. Let \(n\in\N\) and \(z\in\Cs\). In view of \eqref{min0} and \eqref{min1}, we have
 \[
  \sttrhla{\smino{2n-2}{1}}{z}
  =\sttrhla{\smino{2n-1}{1}}{z}
  =-\ek*{\spkuoa{n-1}{1}{z}}\ek*{z\opkuoa{n-1}{1}{z}}^\inv.
 \]
 Using \rtheo{T0820} we obtain furthermore
 \[\begin{split}
  -\ek*{\spkuoa{n-1}{1}{z}}\ek*{z\opkuoa{n-1}{1}{z}}^\inv
  &=-\ek*{z\sphua{n}{z}-\su{0}\ophua{n}{z}}\ek*{z\su{0}^\inv\sphua{n}{z}}^\inv\\
  &=-\su{0}+\frac{1}{z}\su{0}\ek*{\ophua{n}{z}}\ek*{\sphua{n}{z}}^\inv\su{0}.
 \end{split}\]
 Taking additionally into account \eqref{max0}, \eqref{P1058.B1} follows.

In view of \eqref{max0} and \eqref{max1}, we have
 \[
  \sttrhla{\smaxo{2n-1}{1}}{z}
  =\sttrhla{\smaxo{2n}{1}}{z}
  =-\ek*{\sphuoa{n}{1}{z}}\ek*{\ophuoa{n}{1}{z}}^\inv.
 \]
 Using \rtheo{T0820} we obtain furthermore
 \[\begin{split}
  -\ek*{\sphuoa{n}{1}{z}}\ek*{\ophuoa{n}{1}{z}}^\inv
  &=-\ek*{\spkua{n}{z}-\su{0}\opkua{n}{z}}\ek*{\su{0}^\inv\spkua{n}{z}}^\inv\\
  &=-\su{0}+\su{0}\ek*{\opkua{n}{z}}\ek*{\spkua{n}{z}}^\inv\su{0}.
 \end{split}\]
 Taking additionally into account \eqref{min0}, \eqref{P1058.B2} follows.
\eproof

 From the identities derived in the proof of \rprop{P1058}, we can easily obtain the following relations between the matrix polynomials \(\PQPQ\) and \([\seq{\ophuo{k}{1}}{k}{0}{\infi},\seq{\sphuo{k}{1}}{k}{0}{\infi},\seq{\opkuo{k}{1}}{k}{0}{\infi},\seq{\spkuo{k}{1}}{k}{0}{\infi}]\):
\btheol{Tnew1}
 Let \(\seqsinf\in\Kgqinf\) with first \tKt{} $(s_j^\Sta{1})_{j=0}^\infi$ and \tStiqosoomp{} \(\PQPQ\). Then $(s_j^\Sta{1})_{j=0}^\infi$ belongs to \(\Kgqinf\). Denote by \([\seq{\ophuo{k}{1}}{k}{0}{\infi},\seq{\sphuo{k}{1}}{k}{0}{\infi},\seq{\opkuo{k}{1}}{k}{0}{\infi},\seq{\spkuo{k}{1}}{k}{0}{\infi}]\) the \tStiqosoompo{$(s_j^\Sta{1})_{j=0}^\infi$}. Then
\begin{align*}
  \su{0}\ek*{\ophua{n}{z}}\ek*{\sphua{n}{z}}^\inv\su{0}+\ek*{\spkuoa{n-1}{1}{z}}\ek*{\opkuoa{n-1}{1}{z}}^\inv
  &=z\su{0}
 \shortintertext{and}
 \su{0}\ek*{\opkua{n}{z}}\ek*{\spkua{n}{z}}^\inv\su{0}+\ek*{\sphuoa{n}{1}{z}}\ek*{\ophuoa{n}{1}{z}}^\inv
 &=\su{0}
\end{align*}
 for all \(n\in\N\) and all \(z\in\Cs\).
\etheo

 Note that similar interrelations as exposed in \rtheoss{T0820}{Tnew1} between the polynomials \(\PQPQ\) and \([\seq{\ophuo{k}{1}}{k}{0}{\infi},\seq{\sphuo{k}{1}}{k}{0}{\infi},\seq{\opkuo{k}{1}}{k}{0}{\infi},\seq{\spkuo{k}{1}}{k}{0}{\infi}]\) were in the scalar case considered in~\cite{CR16}. \rprop{P1058} can also be seen from the following matrix continued fraction expansions, which appear in connection with matrix Hurwitz type polynomials in~\zita{MR3327132}. For $A,B\in \Cqq $ with $B$ invertible, set \symb{\frac{A}{B}\defg AB^\inv}.

\bpropnl{\zitaa{MR3324594}{\ctheo{3.4}}}{prop2000}
 Let \(\seqsinf\in\Kgqinf\) with \tdsp{} \([\seq{\DSpl{k}}{k}{0}{\infty},\seq{\DSpm{k}}{k}{0}{\infty}]\). For all \(n\in\NO\) and all \(z\in\Cs\), then
 \begin{align*}
 \sttrhla{\smin{2n}}{z}
 &
 =\cfrac{\Iq }{-z\DSpm{0}+
   \cfrac{\Iq }{\DSpl{0}+
    \cfrac{\Iq}{+
     \cfrac{\ddots}{-z\DSpm{n-1}+
      \cfrac{\Iq }{\DSpl{n-1}-z^\inv \DSpm{n}^\inv}}}}}
 \shortintertext{and}
 \sttrhla{\smax{2n}}{z}
  &
  =\cfrac{\Iq }{-z\DSpm{0}+
   \cfrac{\Iq }{\DSpl{0}+
    \cfrac{\Iq }{+
     \cfrac{\ddots}{\DSpl{n-2}+
      \cfrac{\Iq }{-z\DSpm{n-1}+\DSpl{n-1}^\inv}}}}}.
 \end{align*}
\end{prop}


\appendix
\section{Orthogonal matrix polynomials on $\ra$}
 Let us recall some notions on orthogonal matrix polynomials (OMP) which were used in~\zitas{CR13,MR3275438}. Let $P$ be a complex $p \times q$~matrix polynomial. For all $n\in\NO $, let
  \[
    Y^{[P]}_{n}
    \defg
    \begin{bmatrix}
     A_0\\
     A_1\\
     \vdots\\
     A_n
    \end{bmatrix},
\]
 where $(A_j)_{j=0}^\infi$ is the unique sequence of complex $p\times q$~matrices such that for all $z\in\C$ the polynomial $P$ admits the representation $P(z)=\sum_{j=0}^\infi z^jA_j$. Furthermore, we denote by $\deg  P\defg \sup\setaa{j\in\NO}{A_j\neq\Opq}$ the \noti{degree of $P$}{degree}. Observe that in the case $P(z)=\Opq $ for all $z\in\C$ we have thus $\deg P=-\infi$. If $k\defg \deg  P\geq0$, we refer to $A_k$ as the \noti{leading coefficient of $P$}{leading coefficient}.

\begin{defn}\label{def3.2}
 Let $\kappa\in\NOinf $, and let $\seqs{2\kappa}$ be a sequence of complex \tqqa{matrices}.  A sequence $(P_k)_{k=0}^\kappa$ of complex \tqqa{matrix} polynomials is called a \noti{\tmrosmpa{$\seqs{2\kappa}$}}{monic right orthogonal system of matrix polynomials} if the following three conditions are fulfilled:
 \bAeqi{0}
  \il{MLOS.I} $\deg P_k=k$ for all $k\in \mn{0}{\kappa}$.
  \il{MLOS.II} $P_k$ has the leading coefficient $\Iq$ for all $k\in\mn{0}{\kappa}$.
  \il{MLOS.III} $\rk{Y_n^{[P_j]}}^\ad\Hu{n}Y_n^{[P_k]}=   \Oqq$ for all  $j,k\in\mn{0}{\kappa}$ with $j\neq k$, where $n\defg \max\set{j,k}$.
 \eAeqi
\end{defn}
%
\bremnl{cf.~\zitaa{MR3275438}{\cremp{3.6}{1652}}}{R2312}
 Let $\kappa\in \NOinf$ and let $\seqs{2\kappa}$ be a sequence of complex \tqqa{matrices} such that the block Hankel matrix \(\Hu{n}\) is \tpH{} for all \(n\in\mn{0}{\kappa}\). Denote by $(P_k)_{k=0}^\kappa$ the \tmrosmpa{$\seqs{2\kappa}$}. Let $\sigma$ be a \tnnH{} \tqqa{measure} on a non-empty Borel subset $\Omega$ of \(\R\) satisfying $s_j=\int_\Omega t^j\sigma(\dif t)$ for all \(j\in\mn{0}{2\kappa}\). Then
 \[
  \int_\Omega\ek*{P_j(t)}^\ad\sigma(\dif t)\ek*{P_k(t)}
  =
 \begin{cases}
  \Oqq\incase{j\neq k}\\
  \Lu{n}\incase{j=k}
 \end{cases}
 \]
 for all $j,k\in\mn{0}{\kappa}$.
\end{rem}

\bibliography{string_arxiv}

\def\cprime{$'$}
\begin{thebibliography}{10}

\bibitem{MR2155645}
V.~M. Adamyan and I.~M. Tkachenko.
\newblock Solution of the {S}tieltjes truncated matrix moment problem.
\newblock {\em Opuscula Math.}, 25(1):5--24, 2005.

\bibitem{MR2215856}
V.~M. Adamyan and I.~M. Tkachenko.
\newblock General solution of the {S}tieltjes truncated matrix moment problem.
\newblock In {\em Operator theory and indefinite inner product spaces}, volume
  163 of {\em Oper. Theory Adv. Appl.}, pages 1--22. Birkh\"auser, Basel, 2006.

\bibitem{MR0184042}
N.~I. Akhiezer.
\newblock {\em The classical moment problem and some related questions in
  analysis}.
\newblock Translated by N. Kemmer. Hafner Publishing Co., New York, 1965.

\bibitem{MR0290157}
T.~And{\^o}.
\newblock Truncated moment problems for operators.
\newblock {\em Acta Sci. Math. (Szeged)}, 31:319--334, 1970.

\bibitem{MR975671}
V.~A. Bolotnikov.
\newblock Descriptions of solutions of a degenerate moment problem on the axis
  and the halfaxis.
\newblock {\em Teor. Funktsi\u\i\ Funktsional. Anal. i Prilozhen.},
  (50):25--31, i, 1988.

\bibitem{MR1362524}
V.~A. Bolotnikov.
\newblock Degenerate {S}tieltjes moment problem and associated {$J$}-inner
  polynomials.
\newblock {\em Z. Anal. Anwendungen}, 14(3):441--468, 1995.

\bibitem{MR1433234}
V.~A. Bolotnikov.
\newblock On a general moment problem on the half axis.
\newblock {\em Linear Algebra Appl.}, 255:57--112, 1997.

\bibitem{MR1722780}
V.~A. Bolotnikov and L.~A. Sakhnovich.
\newblock On an operator approach to interpolation problems for {S}tieltjes
  functions.
\newblock {\em Integral Equations Operator Theory}, 35(4):423--470, 1999.

\bibitem{MR1807884}
G.-N. Chen and Y.-J. Hu.
\newblock A unified treatment for the matrix {S}tieltjes moment problem in both
  nondegenerate and degenerate cases.
\newblock {\em J. Math. Anal. Appl.}, 254(1):23--34, 2001.

\bibitem{MR1670527}
G.-N. Chen and X.-Q. Li.
\newblock The {N}evanlinna-{P}ick interpolation problems and power moment
  problems for matrix-valued functions.
\newblock {\em Linear Algebra Appl.}, 288(1-3):123--148, 1999.

\bibitem{CR13}
A.~E. Choque~Rivero.
\newblock The {S}tieltjes matrix $n$-convergent via a new {D}yukarev's
  resolvent matrix representation for the nondegenerate finite matricial
  {S}tieltjes moment problem.
\newblock Submitted to Linear Algebra Appl., July 2013.

\bibitem{MR3324594}
A.~E. Choque~Rivero.
\newblock On {D}yukarev's resolvent matrix for a truncated {S}tieltjes matrix
  moment problem under the view of orthogonal matrix polynomials.
\newblock {\em Linear Algebra Appl.}, 474:44--109, 2015.

\bibitem{MR3327132}
A.~E. Choque~Rivero.
\newblock On matrix {H}urwitz type polynomials and their interrelations to
  {S}tieltjes positive definite sequences and orthogonal matrix polynomials.
\newblock {\em Linear Algebra Appl.}, 476:56--84, 2015.

\bibitem{CR16}
A.~E. Choque~Rivero.
\newblock New interrelations between the coefficients and the {M}arkov
  parameters of the polynomial of odd degree and orthogonal polynomials
  generated by the {R}outh-{M}arkov parameters.
\newblock Submitted to IEEE Trans. Autom. Control, 2016.

\bibitem{MR3275438}
A.~E. Choque~Rivero and C.~M{\"a}dler.
\newblock On {H}ankel positive definite perturbations of {H}ankel positive
  definite sequences and interrelations to orthogonal matrix polynomials.
\newblock {\em Complex Anal. Oper. Theory}, 8(8):1645--1698, 2014.

\bibitem{Dyu81}
{\relax Yu}.~M. Dyukarev.
\newblock The {S}tieltjes matrix moment problem.
\newblock Deposited in VINITI (Moscow) at~22.03.81, No.~2628-81, 1981.
\newblock Manuscript, 37~pp.

\bibitem{MR686076}
{\relax Yu}.~M. Dyukarev.
\newblock Multiplicative and additive {S}tieltjes classes of analytic
  matrix-valued functions and interpolation problems connected with them. {II}.
\newblock {\em Teor. Funktsi\u\i\ Funktsional. Anal. i Prilozhen.},
  (38):40--48, 127, 1982.

\bibitem{MR1699439}
{\relax Yu}.~M. Dyukarev.
\newblock A general scheme for solving interpolation problems in the
  {S}tieltjes class that is based on consistent integral representations of
  pairs of nonnegative operators. {I}.
\newblock {\em Mat. Fiz. Anal. Geom.}, 6(1-2):30--54, 1999.

\bibitem{MR2053150}
{\relax Yu}.~M. Dyukarev.
\newblock Indeterminacy criteria for the {S}tieltjes matrix moment problem.
\newblock {\em Mat. Zametki}, 75(1):71--88, 2004.

\bibitem{MR2735313}
{\relax Yu}.~M. Dyukarev, B.~Fritzsche, B.~Kirstein, and C.~M{\"a}dler.
\newblock On truncated matricial {S}tieltjes type moment problems.
\newblock {\em Complex Anal. Oper. Theory}, 4(4):905--951, 2010.

\bibitem{MR2570113}
{\relax Yu}.~M. Dyukarev, B.~Fritzsche, B.~Kirstein, C.~M{\"a}dler, and H.~C.
  Thiele.
\newblock On distinguished solutions of truncated matricial {H}amburger moment
  problems.
\newblock {\em Complex Anal. Oper. Theory}, 3(4):759--834, 2009.

\bibitem{MR645305}
{\relax Yu}.~M. Dyukarev and V.~{\`E}. Katsnel{\cprime}son.
\newblock Multiplicative and additive {S}tieltjes classes of analytic
  matrix-valued functions and interpolation problems connected with them. {I}.
\newblock {\em Teor. Funktsi\u\i\ Funktsional. Anal. i Prilozhen.},
  (36):13--27, 126, 1981.

\bibitem{MR752057}
{\relax Yu}.~M. Dyukarev and V.~{\`E}. Katsnel{\cprime}son.
\newblock Multiplicative and additive {S}tieltjes classes of analytic
  matrix-valued functions, and interpolation problems connected with them.
  {III}.
\newblock {\em Teor. Funktsi\u\i\ Funktsional. Anal. i Prilozhen.},
  (41):64--70, 1984.

\bibitem{MR975253}
B.~Fritzsche and B.~Kirstein.
\newblock Schwache {K}onvergenz nichtnegativ hermitescher {B}orelma\ss e.
\newblock {\em Wiss. Z. Karl-Marx-Univ. Leipzig Math.-Natur. Reihe},
  37(4):375--398, 1988.

\bibitem{MR3014201}
B.~Fritzsche, B.~Kirstein, and C.~M{\"a}dler.
\newblock On a special parametrization of matricial {$\alpha$}-{S}tieltjes
  one-sided non-negative definite sequences.
\newblock In {\em Interpolation, {S}chur functions and moment problems. {II}},
  volume 226 of {\em Oper. Theory Adv. Appl.}, pages 211--250.
  Birkh\"auser/Springer Basel AG, Basel, 2012.

\bibitem{MR3133464}
B.~Fritzsche, B.~Kirstein, and C.~M{\"a}dler.
\newblock Transformations of matricial {$\alpha$}-{S}tieltjes non-negative
  definite sequences.
\newblock {\em Linear Algebra Appl.}, 439(12):3893--3933, 2013.

\bibitem{142arxiv}
B.~Fritzsche, B.~Kirstein, and C.~M{\"a}dler.
\newblock On matrix-valued {S}tieltjes functions with an emphasis on particular
  subclasses.
\newblock {\tt arXiv:1506.01600 [math.CV]}, June 2015.

\bibitem{114arxiv}
B.~Fritzsche, B.~Kirstein, and C.~M{\"a}dler.
\newblock On a simultaneous approach to the even and odd truncated matricial
  {S}tieltjes moment problem~{I}: {A}n $\alpha$-{S}chur-{S}tieltjes-type
  algorithm for sequences of complex matrices.
\newblock {\tt arXiv:1604.07240 [math.CV]}, Apr. 2016.

\bibitem{141arxiv}
B.~Fritzsche, B.~Kirstein, and C.~M{\"a}dler.
\newblock On a simultaneous approach to the even and odd truncated matricial
  {S}tieltjes moment problem~{II}. {A}n $\alpha$-{S}chur-{S}tieltjes-type
  algorithm for sequences of holomorphic matrix-valued functions.
\newblock {\tt arXiv:1604.07629 [math.CV]}, Apr. 2016.

\bibitem{MR3014197}
B.~Fritzsche, B.~Kirstein, C.~M{\"a}dler, and T.~Schwarz.
\newblock On the concept of invertibility for sequences of complex {$p\times
  q$}-matrices and its application to holomorphic {$p\times q$}-matrix-valued
  functions.
\newblock In {\em Interpolation, {S}chur functions and moment problems. {II}},
  volume 226 of {\em Oper. Theory Adv. Appl.}, pages 9--56.
  Birkh\"auser/Springer Basel AG, Basel, 2012.

\bibitem{MR0114338}
F.~R. Gantmacher and M.~G. Kre{\u\i}n.
\newblock {\em Oszillationsmatrizen, {O}szillationskerne und kleine
  {S}chwingungen mechanischer {S}ysteme}.
\newblock Wissenschaftliche Bearbeitung der deutschen Ausgabe: Alfred St\"ohr.
  Mathematische Lehrb\"ucher und Monographien, I. Abteilung, Bd. V.
  Akademie-Verlag, Berlin, 1960.
\newblock English version: Oscillation matrices and kernels and small
  vibrations of mechanical systems, AMS Chelsea Publishing, Providence, RI,
  revised edition, 2002. Translation based on the 1941 Russian original, Edited
  and with a preface by Alex Eremenko.

\bibitem{MR2038751}
Y.-J. Hu and G.-N. Chen.
\newblock A unified treatment for the matrix {S}tieltjes moment problem.
\newblock {\em Linear Algebra Appl.}, 380:227--239, 2004.

\bibitem{MR0044591}
M.~G. Kre{\u\i}n.
\newblock The ideas of {P}. {L}. \v {C}eby\v sev and {A}. {A}. {M}arkov in the
  theory of limiting values of integrals and their further development.
\newblock {\em Uspehi Matem. Nauk (N.S.)}, 6(4 (44)):3--120, 1951.

\bibitem{MR0233157}
M.~G. Kre{\u\i}n.
\newblock The description of all solutions of the truncated power moment
  problem and some problems of operator theory.
\newblock {\em Mat. Issled.}, 2(vyp. 2):114--132, 1967.

\bibitem{MR0458081}
M.~G. Kre{\u\i}n and A.~A. Nudel{\cprime}man.
\newblock {\em The {M}arkov moment problem and extremal problems}.
\newblock American Mathematical Society, Providence, R.I., 1977.
\newblock Ideas and problems of P. L. {\v{C}}eby{\v{s}}ev and A. A. Markov and
  their further development, Translated from the Russian by D. Louvish,
  Translations of Mathematical Monographs, Vol. 50.

\bibitem{MR1627806}
B.~Simon.
\newblock The classical moment problem as a self-adjoint finite difference
  operator.
\newblock {\em Adv. Math.}, 137(1):82--203, 1998.

\bibitem{MR1508159}
T.-J. Stieltjes.
\newblock Recherches sur les fractions continues.
\newblock {\em Ann. Fac. Sci. Toulouse Sci. Math. Sci. Phys.}, 8(4):J1--J122,
  1894.

\bibitem{MR1508160}
T.-J. Stieltjes.
\newblock Recherches sur les fractions continues [{S}uite et fin].
\newblock {\em Ann. Fac. Sci. Toulouse Sci. Math. Sci. Phys.}, 9(1):A5--A47,
  1895.

\end{thebibliography}
\bibliographystyle{abbrv}

\vfill\noindent
\begin{minipage}{240pt}
 A.~E.~Choque-Rivero\\
 Instituto de F\'isica y Matem\'aticas\\
 Universidad Michoacana de San Nicol\'as de Hidalgo\\
 Ciudad Universitaria\\
 Morelia, Mich.\\
 C.P.~58048\\
 M\'exico\\
 \texttt{abdon@ifm.umich.mx}
\end{minipage}
\hspace{10pt}
\begin{minipage}{160pt}
 C.~M\"adler\\
 Universit\"at Leipzig\\
 Mathematisches Institut\\
 Augustusplatz~10\\
 04109~Leipzig\\
 Germany\\
 \texttt{maedler@math.uni-leipzig.de}
\end{minipage}

\end{document}